\documentclass{article}

\usepackage{arxiv}
\usepackage[T1]{fontenc}    % Use 8-bit T1 fonts
\usepackage[utf8]{inputenc} % For UTF-8 encoding
\usepackage[english]{babel} % Language setup
\usepackage{amsmath}        % For mathematical symbols
\usepackage{amsfonts}       % For mathematical fonts
\usepackage{graphicx}       % For including graphics
\usepackage{float}          % For improved floating object control
\usepackage{booktabs}       % For better table formatting
\usepackage{multirow}       % For multi-row table cells
\usepackage{cite}           % For improved citation handling
\usepackage{url}            % Simple URL typesetting
\usepackage{hyperref}       % For hyperlinks
\usepackage{nicefrac}       % Compact symbols for 1/2, etc.
\usepackage{microtype}      % For microtypography adjustments
\usepackage{indentfirst}    % Indent the first paragraph
\usepackage{xcolor}         % For color customization
\usepackage{paralist}       % For inline lists
\usepackage{fancyhdr}       % For fancy headers and footers
\usepackage{etoolbox}       % For advanced programming features
\usepackage{tikz}           % For drawing diagrams
\usetikzlibrary{matrix, positioning, calc, shapes, arrows.meta} % TikZ libraries
\usepackage{algorithm}      % For writing algorithms
\usepackage{algpseudocode}  % For pseudocode in algorithms
\usepackage[pagewise]{lineno} % Line numbers
\usepackage{subcaption}     % Sub-figures
\usepackage{csquotes}       % For context-sensitive quotation facilities
\usepackage{cleveref}       % For enhanced cross-references

\hypersetup{
    colorlinks=true,
    linkcolor=blue,
    filecolor=magenta,
    urlcolor=cyan,
    pdftitle={Your Title Here},
    pdfpagemode=FullScreen,
}

\usepackage{float}

\title{A multi-mesh approach for accurate computation of multi-target functionals in aerodynamics design}
\author{
  Guanghui Hu \\
   Department of Mathematics, Faculty of Science and Technology \\
   State Key Laboratory of Internet of Things for Smart City
   \\
   University of Macau\\
   Zhuhai UM Science \& Technology Research Institute \\
   Zhuhai, Guangdong Province, China\\
   \texttt{garyhu@um.edu.mo} \\
    \And
  Ruo Li \\
   CAPT, LMAM and School of Mathematical Sciences\\
   Peking University\\
   Chongqing Research Institute of Big Data\\ 
   Peking University, Chongqing 401121, China\\
   \texttt{rli@math.pku.edu.cn} \\
 \And
    Jingfeng Wang \\
   Department of Mathematics\\
   Faculty of Science and Technology\\
   University of Macau\\
   \texttt{shankswang953@gmail.mo} \\
 }

\begin{document}
\maketitle

\begin{abstract}
  Aerodynamic optimal design is crucial for enhancing performance of aircrafts, while calculating multi-target functionals through solving dual equations with arbitrary right-hand sides remains challenging. In this paper, a novel multi-target framework of DWR-based mesh refinement is proposed and analyzed. Theoretically, an extrapolation method is generalized to expand multi-variable functionals, which guarantees the dual equations of different objective functionals can be calculated separately. 
Numerically, an algorithm of calculating multi-target functionals is designed based on the multi-mesh approach, which can help to obtain different dual solutions simultaneously.
One feature of our framework is the algorithm is easy to implement with the help of the hierarchical geometry tree structure and the calculation avoids the Galerkin orthogonality naturally.
The framework takes a balance between different targets even when they are not the same orders of magnitude.
While existing approach uses a linear combination of different components in multi-target functionals for adaptation, it introduces additional coefficients for adjusting.
With each component calculated under a dual-consistent scheme, this multi-mesh framework addresses challenges such as the lift-drag ratio and other kinds of multi-target functionals, ensuring smooth convergence and precise calculations of dual solutions.
\end{abstract}

\section{Introduction}
The aerodynamic optimal design has been playing an important role in a variety of areas\cite{martins2022aerodynamic,vasile2020aerodynamic,lee2020wind}.
There have been several methods in the market for solving the PDEs-constrained optimization problem, including gradient-based methods\cite{wu2022gradient,bouhlel2020scalable}, surrogate-based methods\cite{jiang2020surrogate,gramacy2020surrogates}, evolutionary algorithm\cite{ojstersek2020multi,singh2024hybrid} etc. Among all these methods, a common challenge is efficiently calculating the quantity of interest.

The accurate computation of quantities of interest plays a crucial role in the realm of scientific and engineering applications\cite{ojha2021adaptive,leifsson2016multiobjective}. The fidelity of quantities of interest significantly impacts the aerodynamic efficiency, fuel consumption, and overall performance characteristics of the vehicle. 
However, the precise determination of these quantities is inherently complex and computationally intensive, often necessitating substantial computational resources. 

In addressing these computational challenges, mesh adaptation, especially the adjoint-based mesh adaptation, emerges as a powerful tool, pivotal for enhancing the efficiency and accuracy of numerical simulations. As pointed out by \cite{slotnick2014cfd}, adaptive mesh refinement strategies offer the potential for superior accuracy at reduced cost, but have not seen widespread use. 
Over the past decades, there has been considerable achievement in the development of reliable posteriori error estimates and goal-oriented adaptive mesh refinement method\cite{hartmann2002adaptive,hartmann2007adjoint,hartmann2008multitarget,HARTMANN2015754}. 
Techniques like automation\cite{nemec2014toward}, hp-approach\cite{DOLEJSI2021178,dolejvsi2023anisotropic,dolejsi2022}, multi-precision\cite{liu2022mp}, and machine learning \cite{wang2024towards,chen2021output,chen2020output} have been developed in recent years to enhance the computational efficiency and accuracy. 
Based on the Newton-GMG algorithm\cite{li2005multi,hu2010jcp}, we also conducted the h-adaptive refinement\cite{hu2013adaptive} and implemented the DWR-based error estimation\cite{Meng2022,meng2021fourth} within the steady Euler equations. However, the multi-target functionals still present difficulties where further improvements will be discussed in this paper.

There are following issues to calculate the multi-target functionals. For example, the calculation of multi-target objectives is nontrivial where the dual solutions of these non-smooth and nonlinear functionals can not be guaranteed theoretically. Moreover, the dual consistency, which is an important property in developing stable DWR error estimation towards the goal-oriented mesh adaptivity\cite{hartmann2007adjoint,HARTMANN2015754,wang2023towards}, cannot be preserved trivially within the multi-target objective schemes. One typical example is the calculation of lift-to-drag ratio, which is a classical quantity in the aerodynamic optimal design. There have been many approaches proposed in the literature\cite{arora2014lift,sforza2020direct} for the accurate calculation of such a quantity. However, as shown in \cite{dolejsi2022} and our previous experiments\cite{wang2023towards}, the convergence behavior of lift calculation may not be stable with the increase of element size based on the DWR-based mesh adaptation. Not to mention the lift-drag ratio presents additional challenges for the calculation of dual solutions. Solving the dual equations with arbitrary right-hand sides becomes a research concern in recent years\cite{cary2021cfd,cary2022realizing}. Prior to these, there are some alternative techniques that are useful for dealing such problems.
In \cite{hartmann2008multitarget}, the multi-target functionals are transformed to a linear combination of different components, which is widely used in different scenarios\cite{endtmayer2024posteriori,pardo2010multigoal,endtmayer2019mesh}. However, the target functional generated in a linear combined form of different components brings some drawbacks. i). the coefficients introduce manual interventions that are not known in advance for practical applications. ii). the dual consistency cannot be preserved, since the dual solutions originated from the linear combined form of different components cannot match the dual equations where the dual solutions of the single component are demanded. iii). the dual solutions originated from the linear combined form cannot be adopted for calculating the shape derivatives directly.

In this paper, we propose a novel framework for calculating the multi-target functionals with DWR-based mesh refinement method. Theoretically, based on the technique of extrapolation method developed in \cite{venditti2000adjoint}, a multi-variable functional is expanded. Then the derived formula indicates that the dual solutions of different objectives can be calculated separately. Numerically, a multi-mesh framework is developed to obtain the dual solutions simultaneously. With the data structure developed in AFEPack library\cite{AFEPack}, hierarchical geometry tree, the union of different finite volume spaces is equivalent to the union of different trees. It provides an easily implemented algorithm for multi-mesh calculation discussed in this paper. The finite volume space adopted for calculating different dual solutions is the union of different finite volume spaces corresponding with their objective functionals. Thus, it avoids the Galerkin orthogonality naturally. It should be noted that to realize the vision of CFD in 2030\cite{slotnick2014cfd}, solving the dual equations with arbitrary right-hand sides continues as a research focus in different scenarios\cite{cary2021cfd,cary2022realizing}. In this work, the dual equations share the same left-hand side operator since it is derived from the same governing equation, while the right-hand sides are expressed in different forms which are related to specific objective functionals. With the geometrical multi-grid method, the dual solutions are iterated at the same time. Each component derives the dual equations and dual solutions within a dual-consistent scheme. Benefiting from that, the multi-target functional like the lift-drag ratio can converge smoothly with the increase of element size. Besides, the dual solutions corresponding to different objectives are calculated precisely, which provides the foundation for multi-variable optimization. Without manual intervention, the framework can balance different objectives and generate corresponding meshes for calculation.

Different from the single-mesh method that calculates all the components on a functional space, the multi-mesh framework can solve different targets on different functional spaces. 
The multi-mesh approach is widely used in different fields like  thermoelasticity\cite{solin2010monolithic}, the photonic band structure optimization\cite{wu2018multi}, all-electron density functional theory\cite{kuang2023mul}, etc.
The dual solution correspond with specific target emphasizes those elements which contribute to the computation of specific target. As a result, different dual solutions have different distributions on the computational domain. For example, different target functionals on different airfoils are usually considered in multi-airfoil design. The multi-mesh approach shows a great potential to conduct the multi-target optimization in such context.
A significant advantage of multi-mesh method is all the components in the multi-targets functionals can be calculated precisely, leading to a precise result and robust calculation process of complicated multi-targets functionals.

The rest of this paper is organized as follows. In Section 2, we introduce steady Euler equations and the Newton-GMG method. In Section 3, we have a brief review of the single-mesh version of the multi-target DWR-based mesh adaptation method. In Section 4, we elaborate on the multi-target DWR-based mesh adaptation method and the multi-mesh approach for solving this framework. In order to implement this algorithm correctly, issues such as the dual consistency and Galerkin orthogonality are also discussed in this section. Numerical results have been shown in Section 5 where lift-drag ratio and other kinds of multi-target functionals have been tested.
Further improvements are discussed in Section 6.
 
\section{Steady Euler equations}
In the field of aerodynamic design, there are different optimization objectives. For example, minimizing drag or maximizing lift-drag ratio are important research topics. However, the lift-drag ratio is not easy to calculate. The issues come from the theoretical limitation of the DWR-based mesh adaptation method for dealing with multiple target functionals, singularity comes from the fraction term, and weak regularity comes from the lift calculation itself. Moreover, with the development of computational mathematics, multi-target optimization has become more and more important. Other than the lift-drag ratio, the optimization goal also concerns other kinds of aerodynamic parameters like the momentum coefficient, distribution of the pressure, location of shock waves, and so on. Besides, the optimization objective may consist of different components like the practical application concerns the different target functional on different parts of the body.
In this paper, we are going to establish a framework that can resolve such issues.
Firstly, we begin with the introduction of the steady Euler equations and Newton-GMG solver.
\subsection{Basic notations and finite volume discretization}
 In our prior research \cite{HU2016235, hu2016adjoint, hu2011robust}, a robust solver tailored for the steady-state Euler equations has been successfully developed. For two-dimensional inviscid flows, these equations can be expressed in conservative form as
\begin{equation}
    \nabla \cdot \mathcal{F}(\mathbf{u}) = 0 \quad \text{in} \ \Omega,
\end{equation}
where $\mathbf{u}$ denotes the vector of conservative variables, and $\mathcal{F}(\mathbf{u})$ represents the flux vector, detailed as follows:
\begin{equation}
    \mathbf{u} = \left[ \begin{array}{c} \rho \\ \rho u_x \\ \rho u_y \\ E \end{array} \right],
    \quad \text{and} \quad
    \mathcal{F}(\mathbf{u}) = \left[ \begin{array}{cc}
    \rho u_x & \rho u_y \\
    \rho u_x^2 + p & \rho u_x u_y \\
    \rho u_x u_y & \rho u_y^2 + p \\
    u_x(E + p) & u_y(E + p)
    \end{array} \right],
\end{equation}
with the velocity components $(u_x, u_y)^T$, the density $\rho$, pressure $p$, and total energy $E$. The equation of state for an ideal gas is given by
\begin{equation}
    E = \frac{p}{\gamma - 1} + \frac{1}{2} \rho (u_x^2 + u_y^2),
\end{equation}
assuming a specific heat ratio $\gamma = 1.4$.

To solve these equations numerically, the domain $\Omega \subset \mathbb{R}^2$, bounded by $\Gamma$, is discretized into a series of control volumes or elements, $K_i$, through a shape-regular partitioning, denoted as $\mathcal{K}_h$. Here, $K_i$ intersects $K_j$ at $e_{i,j} = \partial K_i \cap \partial K_j$ with $n_{i,j}$ as the unit outward normal vector on $e_{i,j}$ related to $K_i$.

This discretization leads to a weak formulation of the Euler equations within $\mathcal{V}_{H}$as:
\begin{equation}
    R_{H}(\mathbf{u}) = \sum_{K_i\in\mathcal{V}_H} \left( \int_{K_i} \nabla \cdot \mathcal{F}(\mathbf{u}) \, dx \right) = \sum_{K_i\in\mathcal{V}_H} \left( \oint_{e_{i,j} \in \partial K_i} \mathcal{F}(\mathbf{u}) \cdot n_{i,j} \, ds \right) = 0.
\end{equation}

Introducing a numerical flux function $\mathcal{H}(\mathbf{u}_i, \mathbf{u}_j, n_{i,j})$, the discretized form can be written as:
\begin{equation}
\label{EulerDiscrete}
    \sum_{i} \left( \oint_{e_{i,j} \in \partial K_i} \mathcal{H}(\mathbf{u}_i, \mathbf{u}_j, n_{i,j}) \, ds \right) = 0.
\end{equation}

\subsection{Newton-GMG solver}

To tackle the nonlinear Equation \eqref{EulerDiscrete}, we employ the Newton method for linearization, supplemented by a linear multigrid method for solving the equations as proposed in \cite{li2008multigrid}. Initially, we expand Equation \eqref{EulerDiscrete} using the Taylor series and neglect higher-order terms, simplifying the equations to:

\begin{equation}
    \begin{aligned}
        \label{regularized_equation}
        \alpha \left|\!\left|\sum\limits_{i}\sum\limits_{j}\int_{e_{i,j}\in\partial\mathcal{K}_i}\mathcal{H}(\mathbf{u}_i^{(n)},\mathbf{u}_j^{(n)}, n_{i,j})ds \right|\!\right|_{L_1}\Delta \mathbf{u}_i^{(n)}&+\sum\limits_{i}\sum\limits_{j}\int_{e_{i,j}\in\partial\mathcal{K}_i}\Delta \mathbf{u}_i^{(n)}\frac{\partial\mathcal{H}(\mathbf{u}_i^{(n)},\mathbf{u}_j^{(n)}, n_{i,j})}{\partial \mathbf{u}_i^{(n)}}ds\\
&+\sum\limits_{i}\sum\limits_{j}\int_{e_{i,j}\in\partial\mathcal{K}_i}\Delta \mathbf{u}_j^{(n)}\frac{\partial\mathcal{H}(\mathbf{u}_i^{(n)},\mathbf{u}_j^{(n)}, n_{i,j})}{\partial \mathbf{u}_j^{(n)}}ds\\
&=-\sum\limits_{i}\sum\limits_{j}\int_{e_{i,j}\in\partial\mathcal{K}_i}\mathcal{H}(\mathbf{u}_i^{(n)},\mathbf{u}_j^{(n)}, n_{i,j})ds,
    \end{aligned}
\end{equation}
where ${\partial\mathcal{H}(\cdot,\cdot,\cdot)}\slash{\partial \mathbf{u}}$ represents the Jacobian matrix of the numerical flux, and $\Delta \mathbf{u}_i$ represents the increment of conservative variables for the $i$-th element. Following each Newton iteration, the cell average undergoes an update, $\mathbf{u}_i^{(n+1)} = \mathbf{u}_i^{(n)} + \Delta \mathbf{u}_i^{(n)}$.  

A regularization term is incorporated into the simulation to enhance system stability. Initially, the solution substantially deviates from a steady state, but as iterations advance, it converges toward equilibrium, reducing the regularization term to near zero.

For solving Equation \eqref{regularized_equation}, the geometric multigrid technique is adopted. As demonstrated in \cite{wang2023towards}, applying regularization and geometric multigrid methods to both primal and dual equations significantly enhances solver robustness, ensuring reliable dual equation resolution.
\section{Single mesh framework of multi-target DWR-based adaptation}
\label{singleMeshFrame}

In the research area of the airfoil shape optimal design, lift and
drag are two important quantities that need to be considered.
For example, in the context where a complicated structure consists of different parts that need to be optimized, the target functional should be calculated precisely to correspond with different domains. Besides, the lift-drag ratio, which is also a significant indicator, should be calculated with a robust algorithm to guarantee a stable iteration within the optimization process. The lift and drag values are defined as 
\begin{equation}\label{lift_and_drag}
  \mathcal{J}(\mathbf{u})= \int_{\Gamma}p_{\Gamma}(\mathbf{u})\mathbf{n}\cdot \beta,
\end{equation}
where $p_{\Gamma}$ is the pressure on the boundary and $\beta$ in the above formula is given as 
\begin{equation}\beta=\left\{
  \begin{array}{l}
     (\cos\alpha,\sin\alpha)^T/C_{\infty},\text{ for drag calculation}, \\
     (-\sin\alpha,\cos\alpha)^T/C_{\infty},\text{ for lift calculation}.
  \end{array} 
  \right.
\end{equation}
Here $C_{\infty}$ is defined as
$\gamma p_{\infty}Ma_{\infty}^2l/2$, where
$p_{\infty}, Ma_{\infty},l$ denote the far-field pressure, far-field
Mach number and the chord length of the airfoil, respectively. Prior to this work, there have been some techniques to calculate multiple target quantities with single mesh method\cite{hartmann2008multitarget} which we will give a brief review here.

Given $n$ distinct target functionals $F^i(\cdot)$, $i = 1, 2, \cdots, n$, the solutions are calculated on a same mesh which satisfies the tolerance requirement as follows. Find $\mathcal{V}_{h}$, s.t. for $F^i_{h}(u_{h})\in\mathbb{R}$, where $u_h$ is the solution of $R_h(u_h)=0$ on $\mathcal{V}_h$, $|F^i(u_{\infty})-F^i_h(u_h)|<{TOL}_i$. An alternative method is defining a combined target functional $\widetilde{F}(\cdot):=\sum_{i=1}^n \omega_iF^i(\cdot)$. Adopting the single mesh version DWR-based mesh adaptation method with respect to the combined target functional $\widetilde{F}$ until the tolerance requirement is fulfilled for $i= 1,2,\cdots,n$.
From the perspective of the fully discrete method, the target functional can be expanded with Taylor expansion and ignoring the higher order term as
\begin{equation}
    F^i(u_{\infty})=F^i_h(u_h)+\frac{\partial F^i_h}{\partial u_h}(u_{\infty}-u_h).
\end{equation}
Additionally, the residual term can be approximated with the same technique as 
\begin{equation}
    R(u_{\infty}) = R_h(u_h) + \frac{\partial R_h}{\partial u_h}(u_{\infty}-u_h)=0.
\end{equation}

Defining $e_h$ as $u_{\infty}-u_h$, the error estimate of target functional can be approximated with $({\partial F^i_h}/{\partial u_h})e_h$, where $e_h$ can be obtained by solving the equation
\begin{equation}
    \frac{\partial R_h}{\partial u_h}e_h=-R_h(u_h).
\end{equation}

Even though $e_h$ can derive the error of computed target functionals $F^i(\cdot)$, a drawback of this estimation method is that the necessary local information for guiding adaptive mesh refinement hasn't been included. To design error indicators for guiding the target functional $F^i(\cdot)$, $n$ distinctive dual solutions should be obtained. 
It should be noted that for the dual solutions, the left parts are always ${\partial R_h}/{\partial u_h}$, where the right-hand sides are specific target functionals. If the dual equations are solved at the same functional space, the iteration of dual equations can be processed simultaneously until the precision is reached.

In \cite{hartmann2008multitarget}, the coefficient is set as $|\omega_{l}|=10\cdot|\omega_d|$ , where $\omega_l$ is the coefficient of lift and $\omega_d$ is the coefficient of drag respectively. It is demonstrated in the experiment that lift, drag and other variables can be calculated with expected precision. However, multi-target adapted mesh will compromise between the single-target adapted meshes. Then it cannot be as accurate for the individual target functionals as the respective single-target adapted meshes.
During the shape optimization task, shape derivatives play an important role in indicating the shape variation vector and guiding the deformation process, where the accuracy of different dual solutions is demanded to calculate corresponding values.
In order to make sure the accuracy of different dual solutions can be preserved. The multi-mesh framework is then introduced as follows.

\section{A multi-mesh approach for the DWR-based mesh adaptation}

\subsection{A novel multi-target DWR-based mesh adaptation framework}
In this section, we are going to establish the multi-mesh framework for the DWR-based mesh adaptation method. Based on the analysis, the following multi-target functionals can be process with the multi-mesh approach. 
\begin{itemize}
    \item the multi-variable functionals where variables are separable.
 \item linear combination of single variables.
 \item  the composite form of the above two.
    \end{itemize}

    \begin{figure}[h]\centering
        \includegraphics[width=1.0\textwidth,height=0.28\textheight]{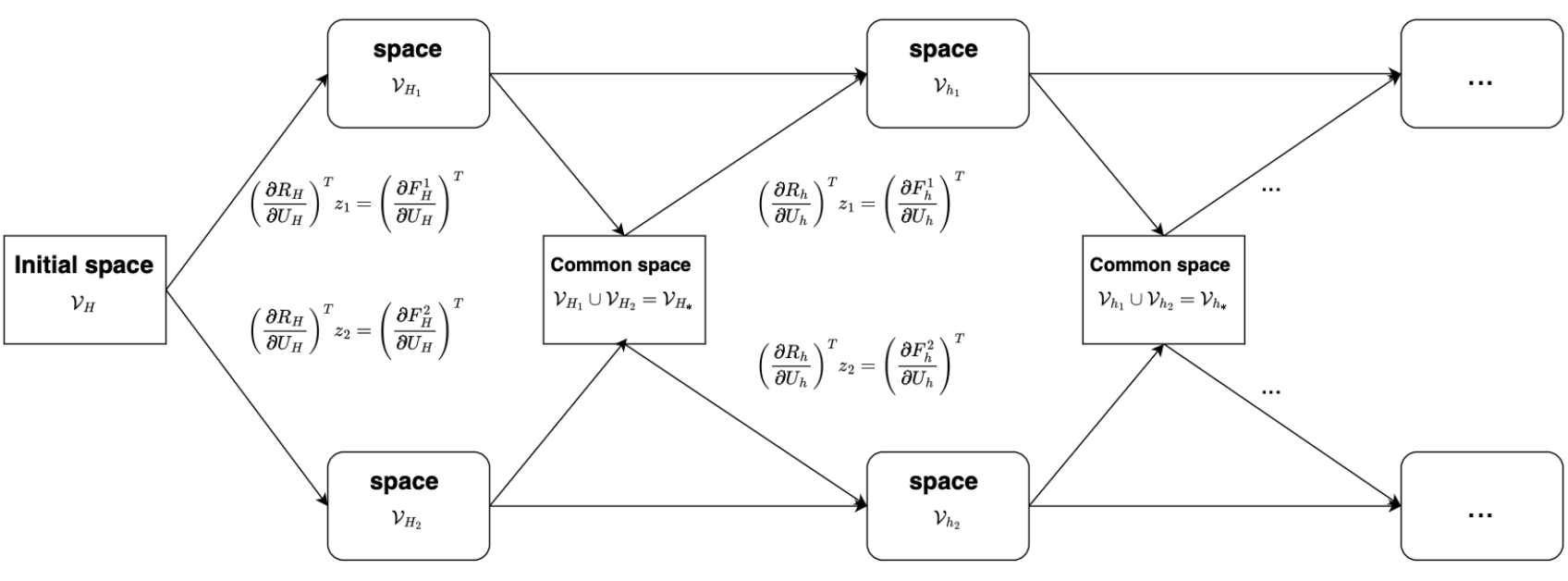}
      \caption{Process of the multi-mesh DWR-based mesh adaptation.}  
      \label{multiMeshProcss}
      \end{figure}

Suppose the target functional $F(x^1,x^2,\cdots,x^n)$ is a multi-variable functional where the variables are separable, i.e.,
\begin{equation}
    F(x^1,x^2,\cdots,x^n)=F^1(x^1)F^2(x^2)\cdots F^n(x^n).
\end{equation}
Denote the distribution of variables on different finite element spaces as $x^1_{H_1}\in \mathcal{V}^1_{H_1}, x^2_{H_2}\in \mathcal{V}^2_{H_2}, \cdots, x^n_{H_n}\in \mathcal{V}^n_{H_n}$. $\mathcal{V}_{H_i}$ represents the functional space with mesh size $H_i$. Similarly, we use $\mathcal{V}_{h_i}$ with mesh size $h_i$ to denote the functional space which is uniformly refined based on $H_i$.  
 $F_{v}(x^1,x^2,\cdots,x^n)$ is used to denote the functional calculated based on the different variables with respect to its own space, i.e.,
\begin{equation}
    F_{v}(x^1_{h_1},x^2_{h_2},\cdots,x^n_{h_n})=F^1_{h_1}(x^1_{h_1})F^2_{h_2}(x^2_{h_2})\cdots F^n_{h_n}(x^n_{h_n}).
\end{equation}
Common functional spaces $\mathcal{V}_{H_*}$ and $\mathcal{V}_{h_*}$ can be obtained by calculating the union of all the functional spaces respectively, i.e.,
\begin{equation}
    \mathcal{V}_{H_*}=\mathcal{V}^1_{H_1}\cup\mathcal{V}^2_{H_2}\cup\cdots\mathcal{V}^n_{H_n}, \qquad \mathcal{V}_{h_*}=\mathcal{V}^1_{h_1}\cup\mathcal{V}^2_{h_2}\cup\cdots\mathcal{V}^n_{h_n}.
\end{equation}
With the multiple variables' Taylor expansion in the first order and ignoring the higher-order term, the following equality can be derived,
\begin{equation}
\label{orignTaylor}
\begin{aligned}
     F_{v}(x^1_{h_1},x^2_{h_2},\cdots,x^n_{h_n})&=F_{v}(I_{h_1}x^1_{H_*}, I_{h_2}x^2_{H_*}, \cdots,I_{h_n}x^n_{H_*})\\
     &+\frac{\partial F_v}{\partial x^1_{h_1}}(I_{h_1}x^1_{H_*}, I_{h_2}x^2_{H_*}, \cdots,I_{h_n}x^n_{H_*})(x^1_{h_1}-I_{h_1}x^1_{H_*}) + \cdots \\
     &+ \frac{\partial F_v}{\partial x^n_{h_n}}(I_{h_1}x^1_{H_*}, I_{h_2}x^2_{H_*}, \cdots,I_{h_n}x^n_{H_*})(x^n_{h_n}-I_{h_n}x^n_{H_*}),
     \end{aligned}
\end{equation}
where $I_{h_i}$ is a projection operator maps variable from $\mathcal{V}_{H_*}$ to $\mathcal{V}_{h_i}$.
Since the variables are separable, then
\begin{equation}
    \frac{\partial F_v}{\partial x^i_{h_i}}(I_{h_1}x^1_{H_*}, I_{h_2}x^2_{H_*}, \cdots,I_{h_n}x^n_{H_*}) = C_i \frac{\partial F^i_{h_i}}{\partial x^i_{h_i}}(I_{h_i}x^i_{H_*}).
\end{equation}
Equation \eqref{orignTaylor} can be converted to
\begin{equation}
\label{convertTaylor}
\begin{aligned}
     F_{v}(x^1_{h_1},x^2_{h_2},\cdots,x^n_{h_n})&=F_{v}(I_{h_1}x^1_{H_*}, I_{h_2}x^2_{H_*}, \cdots,I_{h_n}x^n_{H_*})\\
     &+\sum_{i=1}^n C_i\frac{\partial F^i_{h_i}}{\partial x^i_{h_i}}(I_{h_1}x^1_{H_*}, I_{h_2}x^2_{H_*}, \cdots,I_{h_n}x^n_{H_*})(x^i_{h_i}-I_{h_i}x^j_{H_*}).
     \end{aligned}
\end{equation}
Suppose $R^i_{H_i}(\cdot)$ and $R^i_{h_i}(\cdot)$ are the residual forms of the governing equations on space $\mathcal{V}^i_{H_i}$ and $\mathcal{V}^i_{h_i}$ respectively. The residual of refined space satisfies
\begin{equation}
\label{resG}
    R^i_{h_i}(u^i_{h_i})=0,
\end{equation}
where $u^i_{h_i}$ is the solution of governing equations on space $\mathcal{V}^i_{h_i}$. Specifically, the equation \eqref{convertTaylor} can be expanded from the point of the solution, i.e.,
\begin{equation}
\label{pointTaylor}
\begin{aligned}
     F_{v}(u^1_{h_1},u^2_{h_2},\cdots,u^n_{h_n})&=F_{v}(I_{h_1}u^1_{H_*}, I_{h_2}u^2_{H_*}, \cdots,I_{h_n}u^n_{H_*})\\
     &+\sum_{i=1}^n C_i \frac{\partial F^i_{h_i}}{\partial u^i_{h_i}}(I_{h_1}u^1_{H_*}, I_{h_2}u^2_{H_*}, \cdots,I_{h_n}u^n_{H_*})(u^i_{h_i}-I_{h_i}u^j_{H_*}).
     \end{aligned}
\end{equation}

Similarly, the equation \eqref{resG} can be expanded with Taylor expansion and ignoring the higher order term as
\begin{equation}
\label{resG}
    R^i_{h_i}(u^i_{h_i})=R^i_{h_i}(I_{h_i}u^i_{H_*}) + \frac{\partial R^i_{h_i}}{\partial u^i_{h_i}}(u^i_{h_i}-I_{h_i}u^i_{H_*}) =0.
\end{equation}
Then the equation can be formulated as 
\begin{equation}
\label{finalTaylor}
     F_{v}(u^1_{h_1},u^2_{h_2},\cdots,u^n_{h_n})=F_{v}(I_{h_1}u^1_{H_*}, I_{h_2}u^2_{H_*}, \cdots,I_{h_n}u^n_{H_*})+\sum_{i=1}^n C_i \left( z^i_{h_i} \right)^T R^i_{h_i}(I_{h_i}u^i_{H_*}),
\end{equation}
where $z^i_{h_i}$ satisfies
\begin{equation}
    \left(\frac{\partial R^i_{h_i}}{\partial u^i_{h_i}}\right)^T z^i_{h_i} + \left(\frac{\partial F^i_{h_i}}{\partial u^i_{h_i}}(I_{h_1}u^1_{H_*}, I_{h_2}u^2_{H_*}, \cdots,I_{h_n}u^n_{H_*})\right)^T = 0.
\end{equation}

The second term  $C_i \left( z^i_{h_i} \right)^T R^i_{h_i}(I_{h_i}u^i_{H_*})$ can be adopted as an indicator to refine mesh concerning the functional $F^i(x^i)$. The coefficient $C_i$ is fixed when concerning the single mesh refinement concerning $H_i$. 
In \cite{wang2024towards}, we established an automatic method to choose a suitable tolerance for the mesh adaptation. The algorithm can deal with indicators with different orders of magnitude. Thus, even if the coefficient $C_i$ is unknown, the mesh adaptation process will not be affected. Taking the target functional with two components as an example, the process can be explained in Figure \ref{multiMeshProcss} and Algorithm \ref{alg:multi-mesh-dwr}, where the index $(k)$ denotes the refined times.

\begin{algorithm}
    \caption{Multi-mesh DWR-based Mesh Adaptation Process}
    \begin{algorithmic}[1]
    
    \Statex \hspace*{\dimexpr-\algorithmicindent}\textbf{Input:} Initial mesh with space $\mathcal{V}_{H}$, multi-variable target functional $F(x^1, x^2, \cdots, x^n)$, tolerance $TOL$
    \Statex \hspace*{\dimexpr-\algorithmicindent}\textbf{Output:} target functional on refined mesh $F_{H^{(k+1)}_*}(x^1, x^2, \cdots, x^n)$
    
    \For{each variable $x^i$ in $x^1, x^2, \cdots, x^n$}
        \State Perform single mesh DWR refinement for $x_i$ to obtain space $\mathcal{V}_{H^{(0)}_i}$
    \EndFor
    
    \Repeat
        \State Merge spaces $\mathcal{V}_{H^{(k)}_i}$ to obtain combined space $\mathcal{V}_{H^{(k)}_*}$
        \State On $\mathcal{V}_{H^{(k)}_*}$, iterate with GMG algorithm to compute dual solutions $z^i_{H_{*}^{(k)}}$ for each variable
        \For{each variable $x_i$ in $x^1, x^2, \cdots, x^n$}
            \State Interpolate $z^i_{H_{*}^{(k)}}$ back to the original space $\mathcal{V}_{H^{(k)}_i}$
            \State Calculate the indicator with residual $R^i_{H_i^{(k)}}(u^i_{H_i^{(k)}})=0$ and $z^i_{H_{*}^{k}}$
            \State Refine the mesh with respect to $x^i$ based on the indicator to obtain space $\mathcal{V}_{H_i^{(k+1)}}$
        \EndFor
    \Until{$||F_{H^{(k)}_*}(x^1, x^2, \cdots, x^n)- F_{H^{(k+1)}_*}(x^1, x^2, \cdots, x^n)|| < TOL$ }
    
    \end{algorithmic}
    \label{alg:multi-mesh-dwr}
\end{algorithm}

\subsection{Hierarchical Geometry Tree}
The mesh refinement algorithm in this work is conducted with the library AFEPack. In \cite{li2005multi}, an efficient and robust management of mesh grids is established with the tree structure in AFEPack. 
The implementation of a hierarchical geometry tree provides several critical advantages. The hierarchical representation allows for detailed and flexible access to geometric information across different dimensions. This capability is particularly beneficial for applications requiring frequent and dynamic adjustments to the mesh, such as adaptive mesh refinement. Moreover, the tree-based storage framework ensures that operations like refinement and coarsening are computationally efficient and systematically organized. Consequently, the hierarchical geometry tree proves to be a powerful tool in the management of mesh grids, supporting complex computational simulations and algorithms with high efficiency and flexibility\cite{kuang2023mul,bao2012h, hu2009multi,di2009computation}. 

\begin{figure}[!h]\centering
\includegraphics[width=0.59\textwidth]{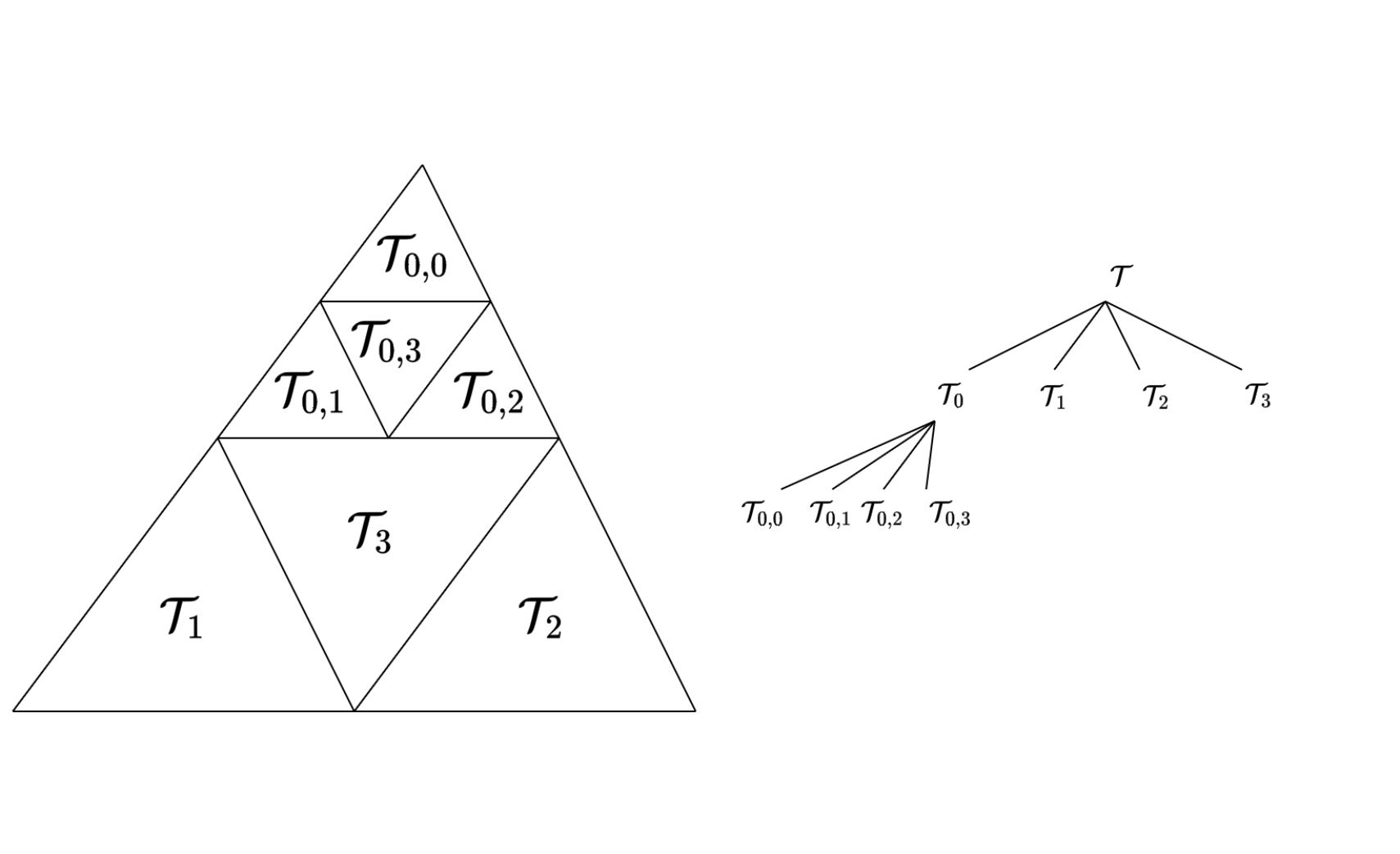}
\includegraphics[width=0.4\textwidth]{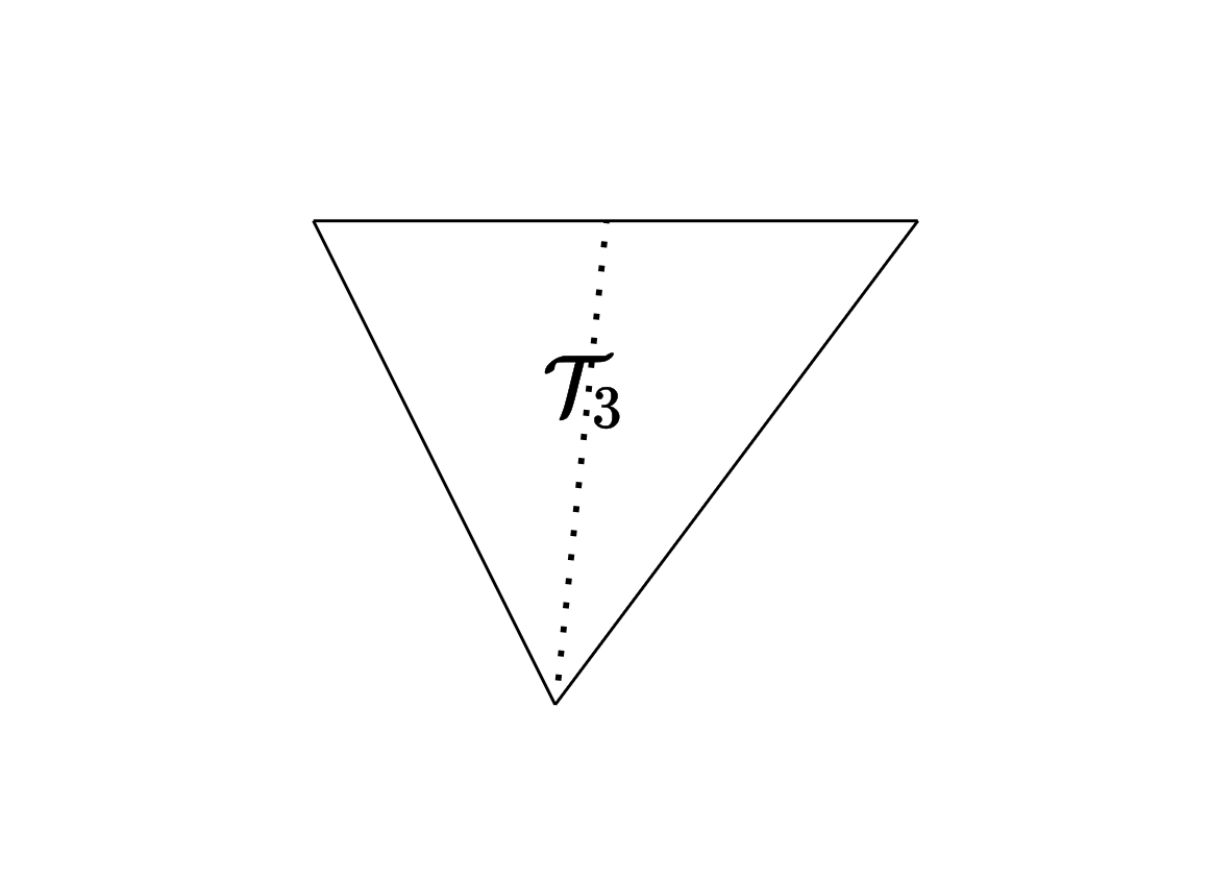}
\caption{Left: hierarchical mesh and hierarchical tree. Right: Twin triangle.}
\label{T1}
\end{figure}

\begin{figure}[!h]\centering
\includegraphics[width=0.33\textwidth]{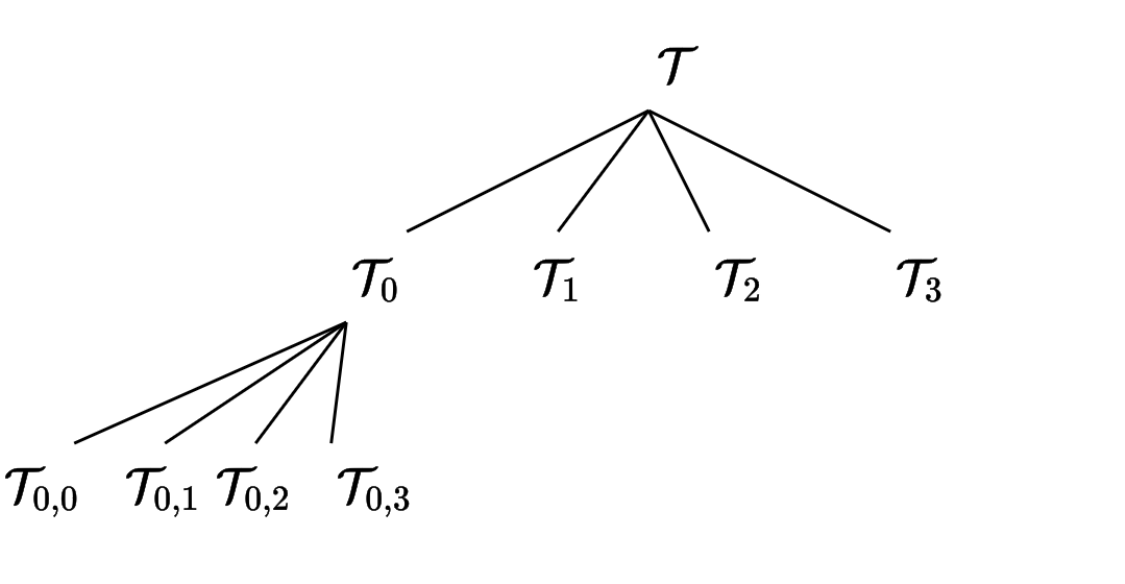}
\includegraphics[width=0.33\textwidth]{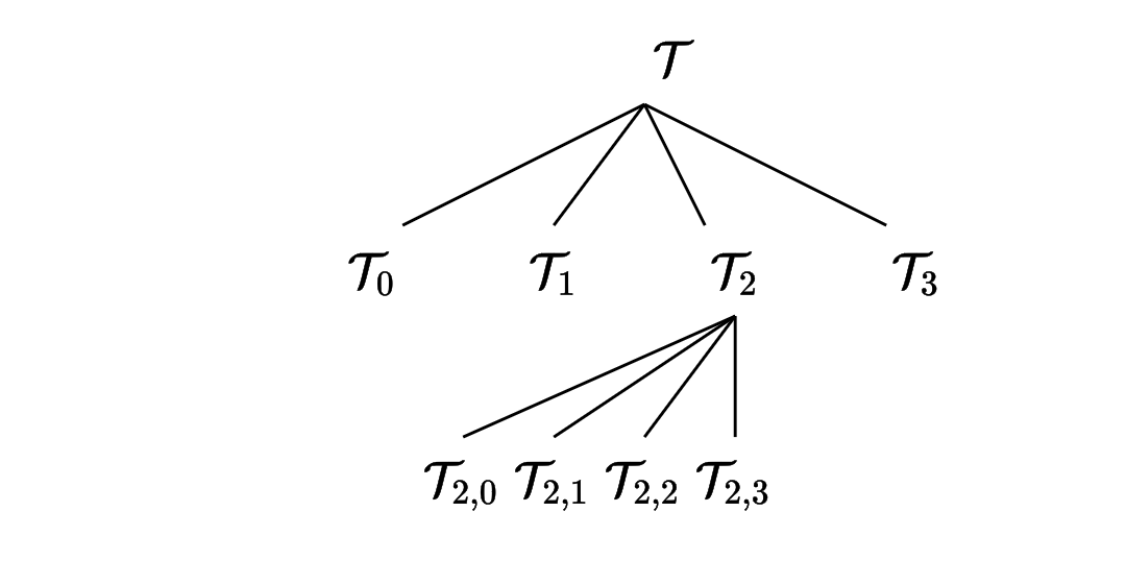}
\includegraphics[width=0.33\textwidth]{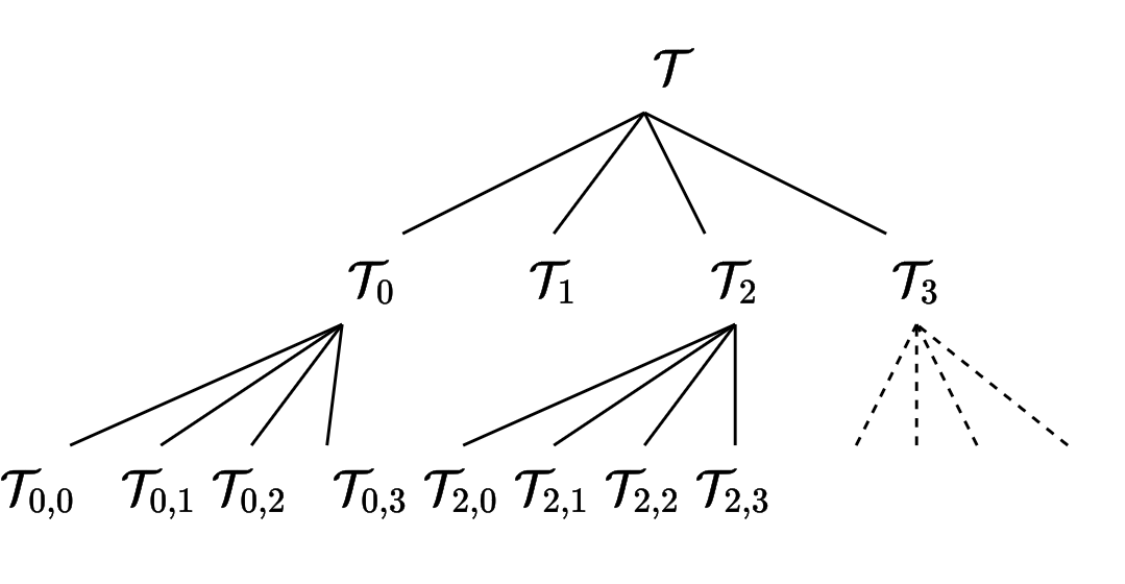}\\
\includegraphics[width=0.33\textwidth]{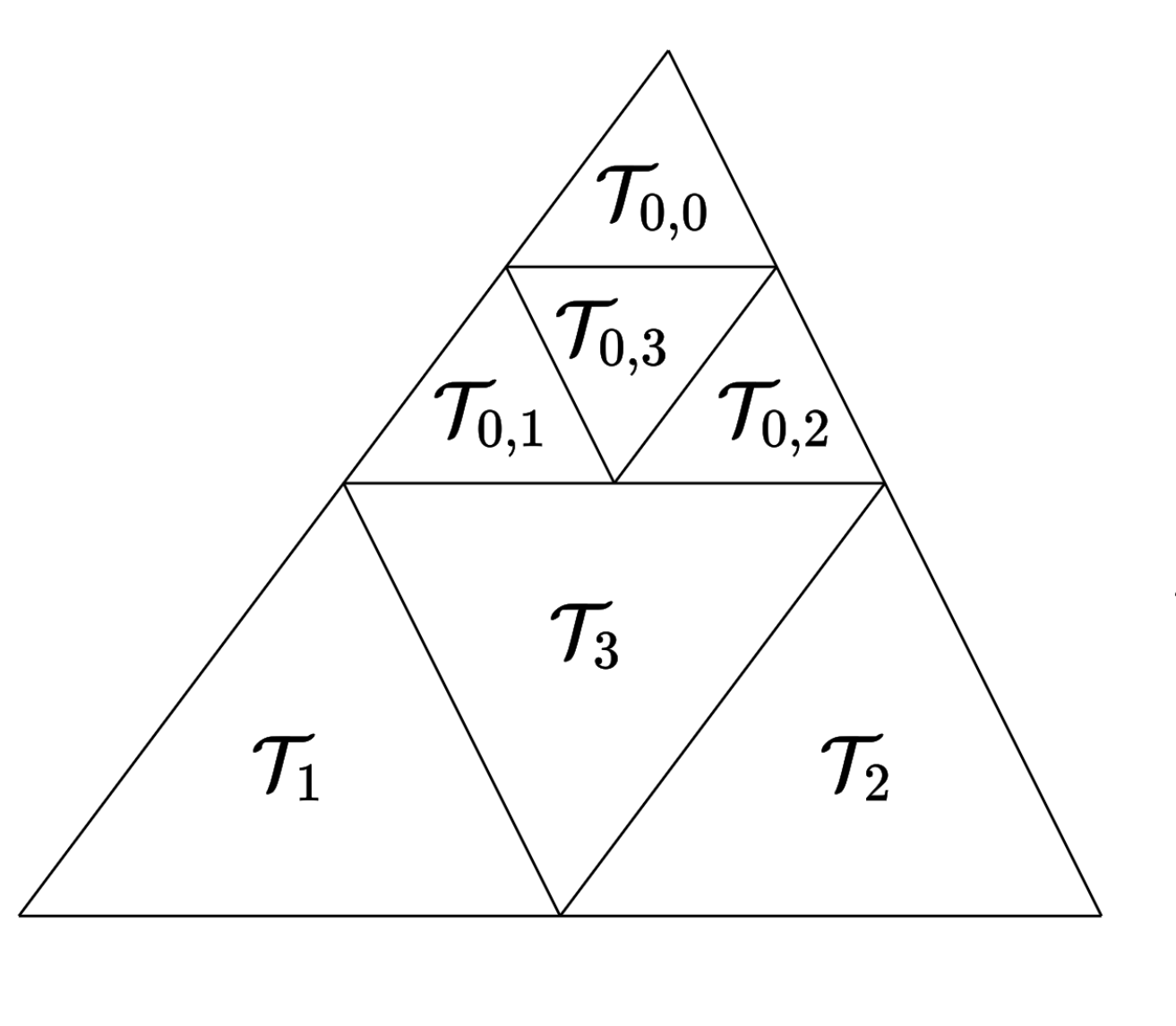}
\includegraphics[width=0.33\textwidth]{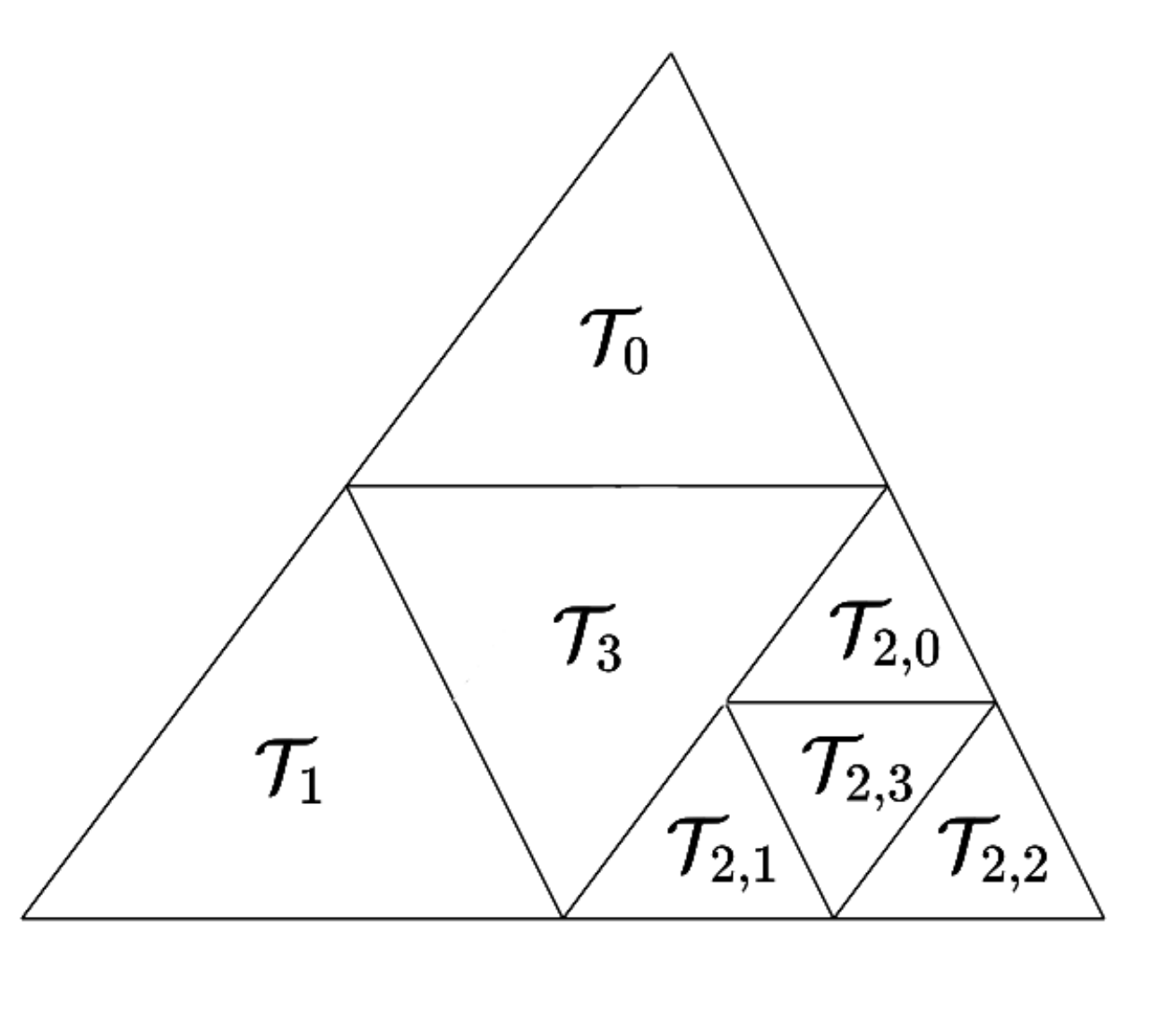}
\includegraphics[width=0.33\textwidth]{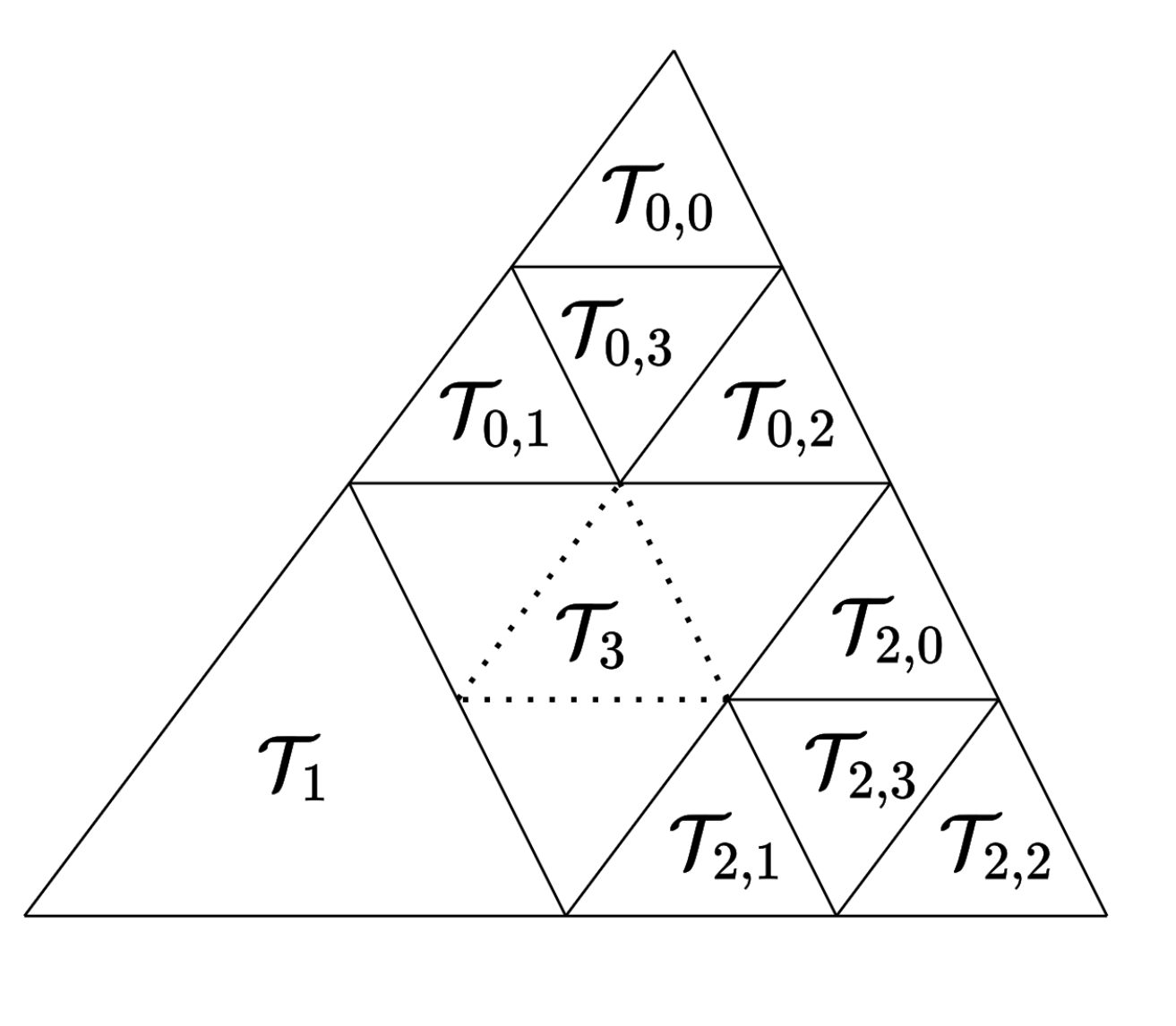}
\caption{The union of tree structure generates the mesh for common space.}
\label{TreeUnion}
\end{figure}

As shown in Fig \ref{T1}, a specific element $\mathcal{T}$ is refined. Meanwhile, the node $\mathcal{T}$ in the tree structure will have $4$ child elements $\mathcal{T}_0,\mathcal{T}_1,\mathcal{T}_2,\mathcal{T}_3$. Subsequently, $\mathcal{T}_0$ gets refined and $4$ child elements $\mathcal{T}_{0,0},\mathcal{T}_{0,1},\mathcal{T}_{0,2},\mathcal{T}_{0,3}$ are obtained. The mesh structure is organized in a hierarchical manner where it begins with zero-dimensional points to two-dimensional triangles. This hierarchical organization enables flexible referencing and manipulation of geometric information, thereby enhancing the efficiency of mesh refinement and coarsening processes. Besides, it should be noted that a hanging point will appear when $\mathcal{T}_0$ is refined. Then the element $\mathcal{T}_3$ will contain four points. In the implementation of the AFEPack library, for two-dimensional mesh refinement, each element is allowed to have at most one hanging point. If more than two hanging points are detected, the element is further refined. The element has one hanging point as shown in the right part of Fig \ref{T1} is named a twin triangle. As illustrated in the previous section, a common space needs to be established to solve the dual solutions with respect to different objectives simultaneously. Within the framework of the finite volume method, this requirement can be met by computing the union of the two meshes. The union of the two meshes is equivalent to merging the tree structures of each element. The union process can seen from Fig \ref{TreeUnion}. As $\mathcal{T}_0$ and $\mathcal{T}_2$ get refined, there appear two hanging points on $\mathcal{T}_3$. Then $\mathcal{T}_3$ should be refined as well.

\subsection{Dual consistency and Galerkin orthogonality}
Concerning the implementation of a multi-mesh DWR-based mesh adaptation method, issues like the dual consistency and Galerkin orthogonality will greatly influence the effect of the algorithm. However, the multi-mesh framework for the DWR-based mesh adaptation can naturally avoid the issues that come from the Galerkin orthogonality. Details will be explained in this section.

Firstly, we start to review the concept of dual consistency, which is closely related to the smoothness of the
discrete dual solutions. If the discretization is implemented under a
dual-consistent scheme, the discrete dual solutions should approximate
the continuous dual solutions as the refinement level
increases. Conversely, a dual-inconsistent scheme may generate dual
solutions with unexpected oscillations or exhibit some nonsmoothness. 

 In
\cite{hartmann2007adjoint}, Hartmann developed analyses about dual
consistency under the discontinuous Galerkin scheme. Motivated by
\cite{hartmann2007adjoint}, further analyses are made to discuss the
dual consistency within Newton-GMG framework in \cite{wang2023towards}.

Suppose the primal solutions are obtained from space $\mathcal{V}^i_{H_{i}}$, the discretized primal equations \eqref{EulerDiscrete} are
denoted as residual form: Find $u^i_{H_{i}} \in\mathcal{V}^i_{H_{i}}$, s.t.
\begin{equation}\label{discretized_primal}
  R_{H_i}^i(u^i_{H_{i}},v^0)=0,\qquad \forall v^0\in\mathcal{V}_{H_i},
\end{equation}
where $v^0$ denotes the test function with $0$ degree.
The notation $R_{H_i}^i(u^i_{H_{i}},v^0)$ is defined as 
\begin{equation}
  \label{discrete_operator}
    R_{H_i}^i(u^i_{H_{i}},v^0)
   :=\displaystyle\sum\limits_{K\in\mathcal{V}_{H_i}}\left\{\int_{\partial K\backslash \Gamma}\mathcal H((u^i_{H_{i}})^+,(u^i_{H_{i}})^-,n_K)(v^0)^++\int_{\partial K\cap\Gamma}\tilde{\mathcal H}((u^i_{H_{i}})^+,\Gamma((u^i_{H_{i}})^+),n_K)(v^0)^+\right\}
\end{equation}

Here in each control volume $K\in\mathcal{K}_h$, we use $(u^i_{H_{i}})^+$ and $(u^i_{H_{i}})^-$ to denote the interior and exterior traces of $u^i_{H_{i}}$ of specific element $K$, respectively. $n_K$ is the corresponding vector from interior to exterior direction.
The continuous version of primal equations
can be defined as: Find $u_{\infty}\in\mathcal{V}$, s.t.
\begin{equation}
\label{primalInf}
  R(u_{\infty},v)=0,\qquad\forall v^0\in\mathcal{V}.
\end{equation}
The notation $R(u_{\infty},v)$ is defined as 
\begin{equation}
    R(u_{\infty},v):=\int_{\Omega} (\nabla\cdot F(u_{\infty}))v^0dx.
\end{equation}

With the  Fr\'{e}chet derivatives, we can denote the continuous exact dual equations of functional $F^i$ as: Find $z_{\infty}\in\mathcal{V}$, s.t.
\begin{equation}
  R'[u_{\infty}](w,z_{\infty})+(F^i)'[u_{\infty}](w)=0,\qquad \forall w\in\mathcal{V}.
\end{equation}
where $z_{\infty}$ is the exact dual solution of the continuous dual equations. Similarly, the discretized dual equations can be defined as :
Find $z_{H_i}\in\mathcal{V}_{H_i}$, s.t.
\begin{equation}
  {R'_{H_i}}[u_{H_i}](w,z_{H_i})+(F^i_{H_i})'[u_{H_i}](w)=0,\qquad \forall w_{H_i}\in\mathcal{V}_{H_i}.
\end{equation}

The primal consistency is held when the exact solution $u_{\infty}$ of equation \eqref{primalInf} 
satisfies the discretized operator:
\begin{equation}
  R_{H_i}^i(u_{\infty},v^0)= 0,\qquad \forall v^0\in\mathcal{V}_0^{H_i}.
\end{equation}

The quantity of interest is defined as \it
dual-consistent\rm\cite{hartmann2007adjoint} with the governing
equations if the discretized operators satisfy:

\begin{equation}\label{DUALConsist}
  R_{H_i}^i[u_{\infty}](w,z)+(F^i_{H_i})'[u_{\infty}](w)=0.
\end{equation}

In \cite{wang2023towards}, the dual consistency of the Newton-GMG framework has been discussed. 
Based on the multi-mesh framework for the DWR-based mesh adaptation, each component of the target functional should be solved within the dual-consistent algorithm.
Other issues come from the error indicators. Considering the error estimate of specific functional $F^i$.
\begin{equation}
\begin{aligned}
    F^i_{H_i}(u^i_{H_i})-F^i_{H_i}(u_\infty)&=\int_0^1 (F^i_{H_i})'[\theta u^i_{H_i} + (1-\theta) u_{\infty}](u^i_{H_i} - u_{\infty}) d\theta\\
    &=-\int_0^1 (R^i_{H_i})'[\theta u^i_{H_i} + (1-\theta) u_{\infty}](u^i_{H_i} - u_{\infty}, z_{\infty}) d\theta \\
    &=-\int_0^1 (R^i_{H_i})'[\theta u^i_{H_i} + (1-\theta) u_{\infty}](u^i_{H_i} - u_{\infty}, z_{\infty}- z^i_{H_i}) d\theta\\
    &=R^i_{H_i}(u^i_{H_i}, z_{\infty})-R^i_{H_i}(u^i_{H_i}, z_{H_{i}})
    \end{aligned}
\end{equation}

It should be noted that if the discretized dual solution $z_{H_i}$ is solved on the same space as $u^i_{H_i}$, the error estimate is identical to zero. Thus, the dual solutions and primal residuals should be calculated on different spaces. According to the algorithm as shown in Algorithm \ref{alg:multi-mesh-dwr}. The dual solutions are calculated on the common spaces, which is the union of different functional spaces corresponding with different target functional. It can avoid the potential issues originating from the Galerkin orthogonality naturally.

\section{Numerical results}
\subsection{Multi-airfoil}
In this section, we are going to test the algorithm with the multi-airfoil case. From the perspective of practical application, the optimization objectives may consist of different parts that need to be calculated precisely. Thus, the multi-airfoil problems approximate this circumstance. Most importantly, the lift and drag have a strong internal connection. Even if the mesh adaptation is conducted towards the lift, the calculation of drag will benefit from the mesh generated from the lift-based adaptation. To mitigate the potential influence from the internal connection, we conduct the multi-target calculation with the multi-airfoil example with the following configurations.
\begin{itemize}
  \item  A domain surrounded by an outer
circle with a radius of 35;
\item  Mach number 0.729, and attack angle
2.31$^\circ$; 
\item Lax-Friedrichs numerical flux;
\item  NACA0012 at $(5,5)$, RAE2822 at $(-5,-5)$;
\item Mirror
  reflection as the solid wall boundary condition;
\item Objectives: Lift coefficient of NACA0012, Drag coefficient of RAE2822.
\end{itemize}

\begin{figure}[!h]\centering
\frame{\includegraphics[width=0.33\textwidth]{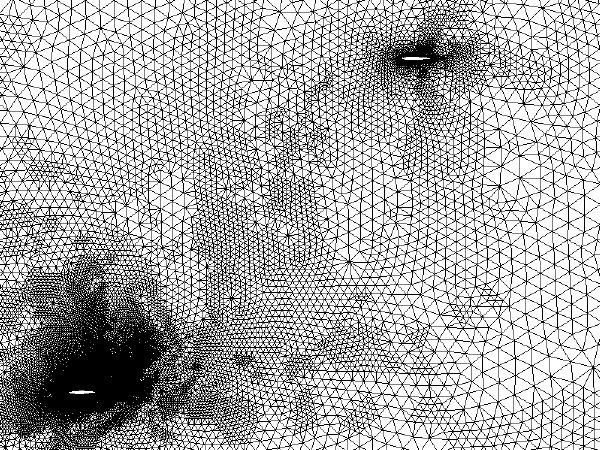}}  
\frame{\includegraphics[width=0.33\textwidth]{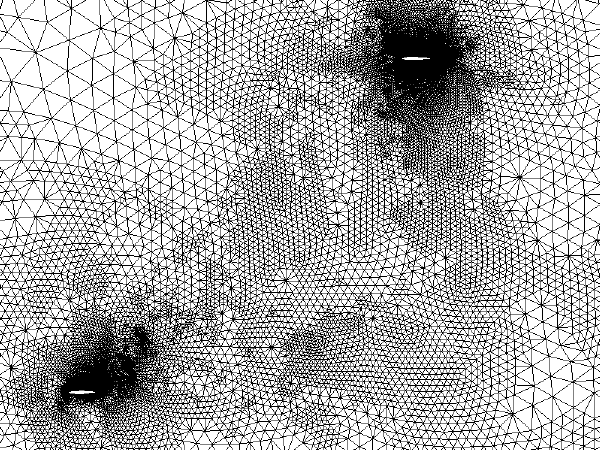}}
\frame{\includegraphics[width=0.33\textwidth]{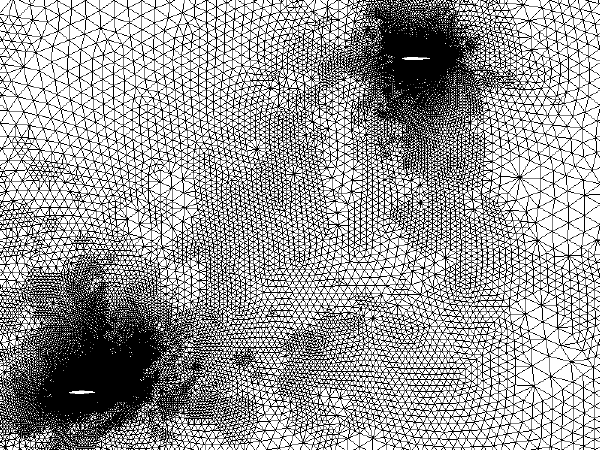}}
\caption{Left: the mesh generated with lift as the target. Middle: the mesh generated with drag as the target. Right: the common mesh.}
\label{multiMeshdwr}
\end{figure}

The process of multi-mesh DWR-based mesh adaptation can be seen in Fig. \ref{multiMeshdwr}. The left part is generated from the DWR-based mesh adaptation with target functional as the drag of RAE2822 at the bottom left while the middle part is generated from the DWR-based mesh adaptation with target functional as the lift of NACA0012 at the top right. The right part is the common mesh for the subsequent dual calculation. 

\begin{figure}[!h]\centering
\frame{\includegraphics[width=0.45\textwidth]{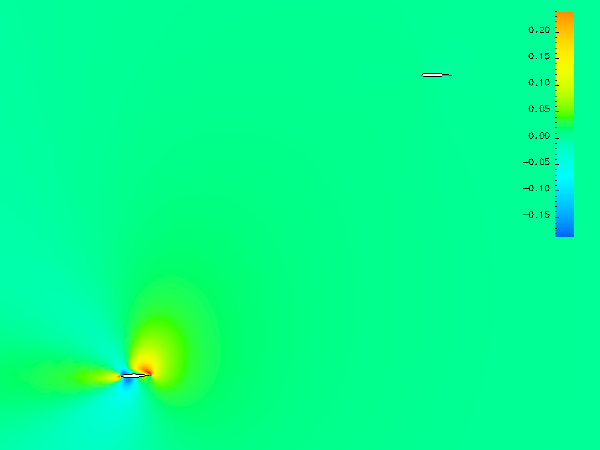}}  
\frame{\includegraphics[width=0.45\textwidth]{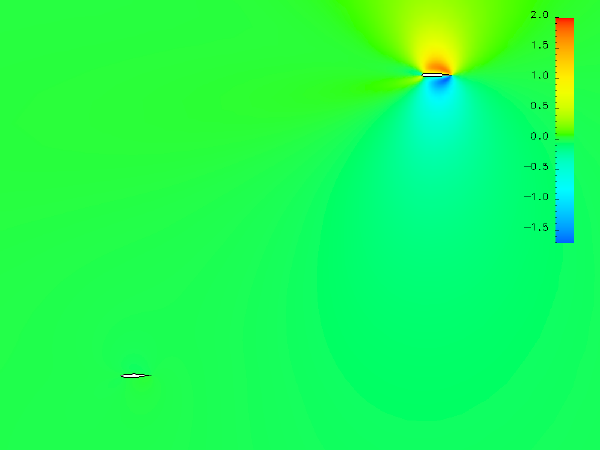}}
\caption{Left: the dual solution of drag. Right: the dual solution of lift.}
\label{multiDual}
\end{figure}

The dual solutions are obtained from the common mesh simultaneously. As shown in Fig. \ref{multiDual}, the dual equations of both target functional are iterated at the same time. The left part is the dual solution of drag while the middle part is that of the lift. 


\begin{figure}[!h]\centering
\includegraphics[width=0.45\textwidth]{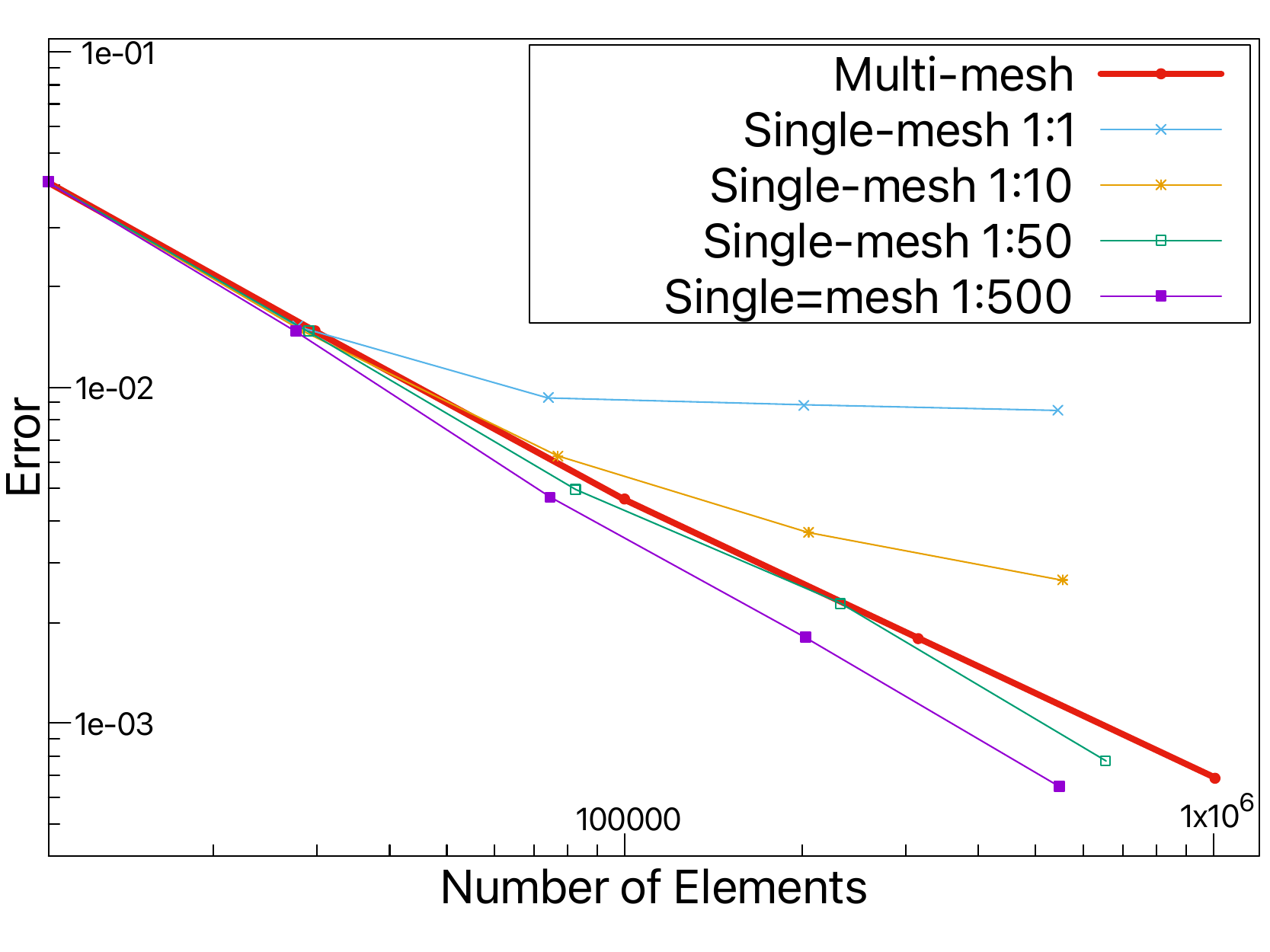}  
\includegraphics[width=0.45\textwidth]{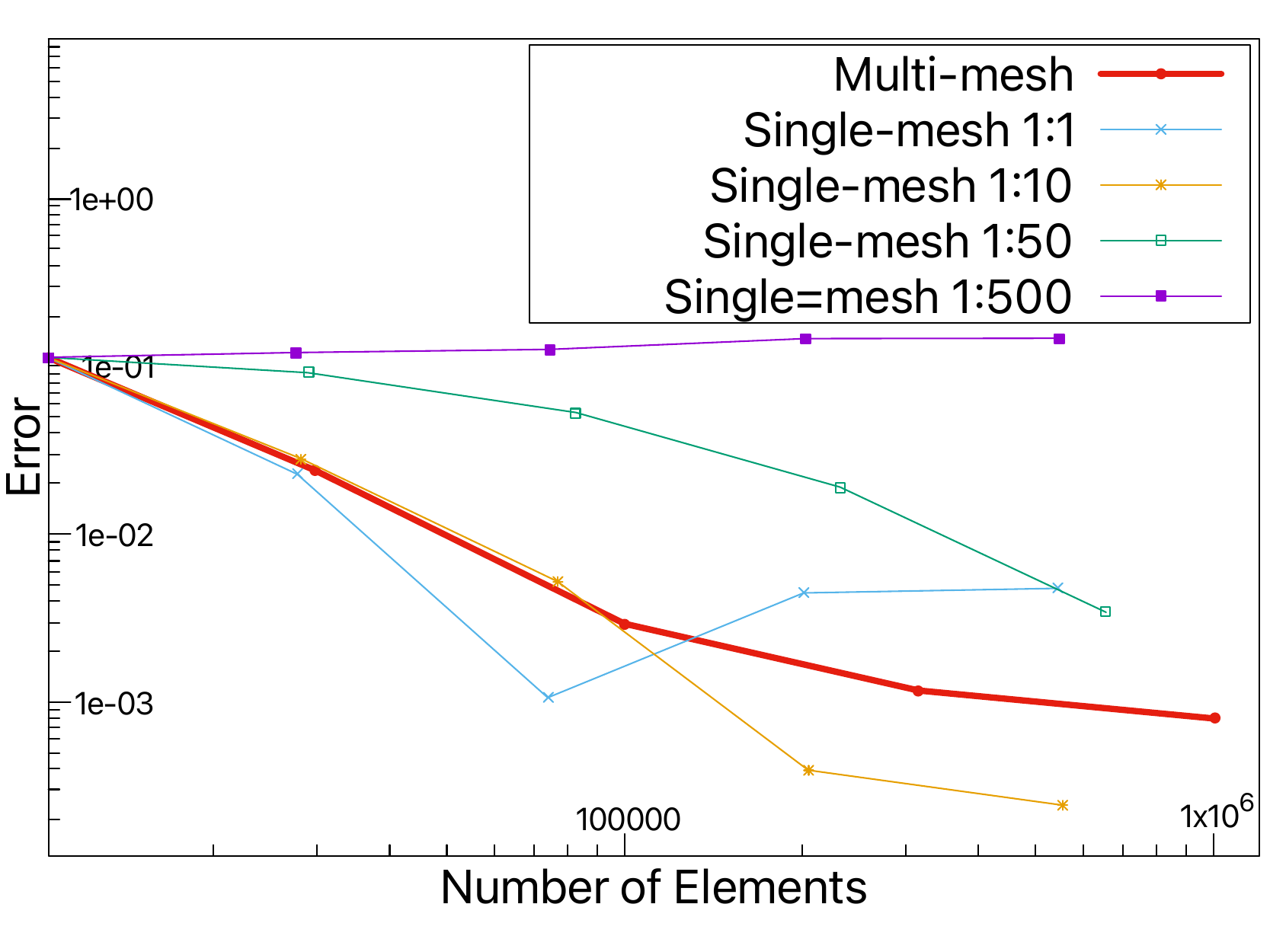}
\caption{Comparison of convergence of multi-mesh method and single mesh method with different coefficients. Left: drag. Right: lift. }
\label{compareSingle}
\end{figure}

\begin{figure}[!h]\centering
    \begin{subfigure}{0.45\textwidth}
        \frame{\includegraphics[width=\textwidth]{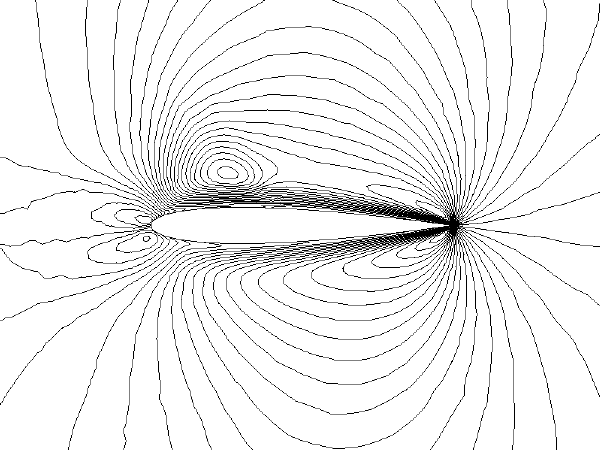}}
        \caption{multi-mesh method}
    \end{subfigure}
    \hfill
    \begin{subfigure}{0.45\textwidth}
        \frame{\includegraphics[width=\textwidth]{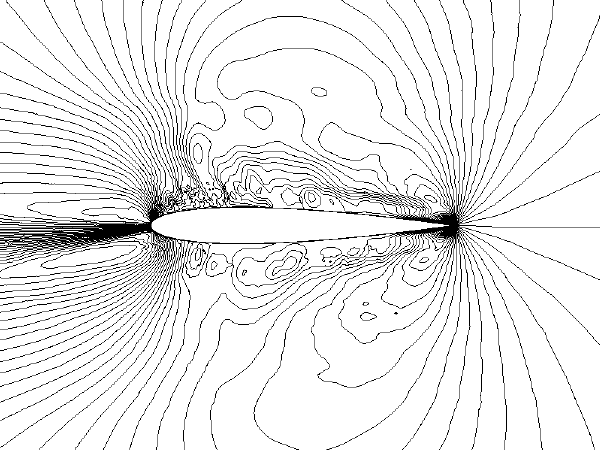}}
        \caption{single mesh method with coefficient $|\omega_l|:|\omega_d| = 1:1$}
    \end{subfigure}
    \\
    \begin{subfigure}{0.45\textwidth}
        \frame{\includegraphics[width=\textwidth]{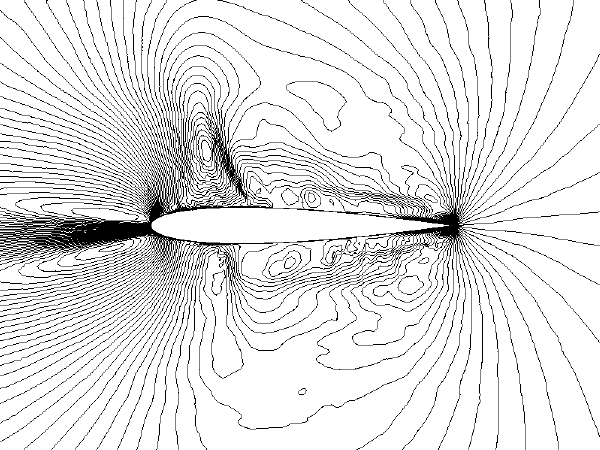}}
        \caption{single mesh method with coefficient  $|\omega_l|:|\omega_d| = 1:10$}
    \end{subfigure}
    \hfill
    \begin{subfigure}{0.45\textwidth}
        \frame{\includegraphics[width=\textwidth]{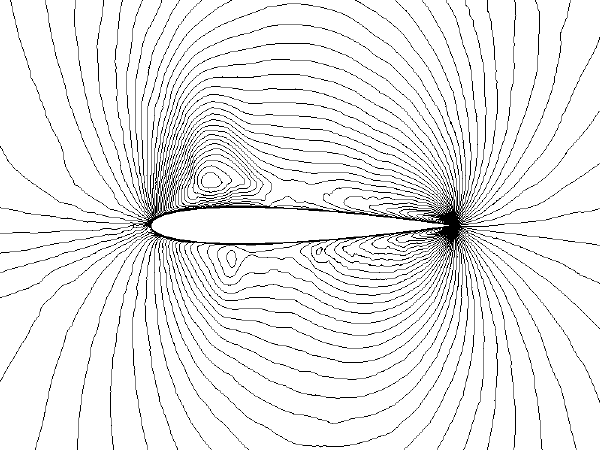}}
        \caption{single mesh method with coefficient  $|\omega_l|:|\omega_d| = 1:50$}
    \end{subfigure}

    \caption{Isolines of the dual solutions around the NACA0012. }
    \label{isonLine}
\end{figure}

The initial mesh gets uniformly refined with $5$ rounds to produce a reference value, which is on a mesh with $10,723,328$ elements where the lift of NACA0012 is $0.564174$ and drag of RAE2822 is $0.0131472$. Then we compare the result of the multi-mesh framework with the single mesh framework where the lift and drag's linear combination $\omega_lF^1+\omega_lF^2$ is adopted for the mesh generation. Here $F^1$ denotes the lift of NACA0012 while $F^2$ denotes the drag of RAE2822. The convergence history of drag and lift, as shown in Fig. \ref{compareSingle}, the multi-mesh method can converge robustly while the behavior of the single mesh method is influenced by the ratio of coefficient significantly. Even though the single method can outperform the multi-mesh one for some cases, the convergence history may not be preserved for other target functional as well. Worse still, it may be oscillatory for a certain ratio. It can be attributed to the issues of dual consistency. The target functional for the dual equations generated from the single mesh method is the combination of lift and drag, not the lift and drag themselves. As shown in Fig.\ref{isonLine}, the isolines of the multi-mesh method and the single mesh method around the NACA0012 airfoil are different. While the multi-mesh method preserves the dual consistency strictly, the isolines are smoother than that of the single-mesh method.
More importantly, as illustrated in previous research experiments\cite{wang2023towards,dolejvsi2023anisotropic}, the convergence history can be more stable with the dual consistency preserved. Even though the convergence rate of the multi-mesh method is smaller in this experiment, which due to the adaptation of two different areas, the precision of the multi-target functional can meet the expectation with mesh adaptation stepped.

\begin{figure}[!h]\centering
    \begin{subfigure}{0.45\textwidth}
        \frame{\includegraphics[width=\textwidth]{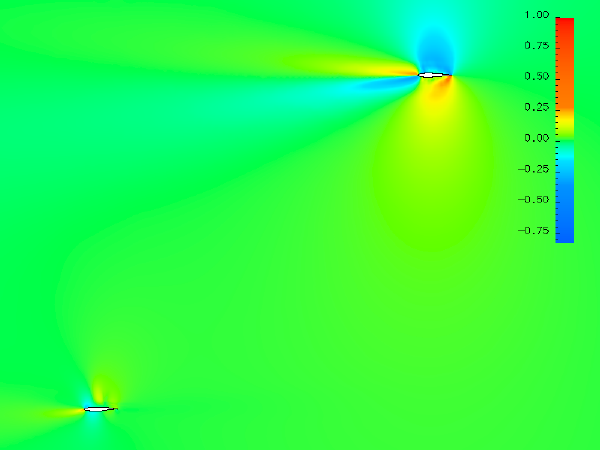}}
        \caption{$|\omega_l|:|\omega_d| = 1:1$}
    \end{subfigure}
    \hfill
    \begin{subfigure}{0.45\textwidth}
        \frame{\includegraphics[width=\textwidth]{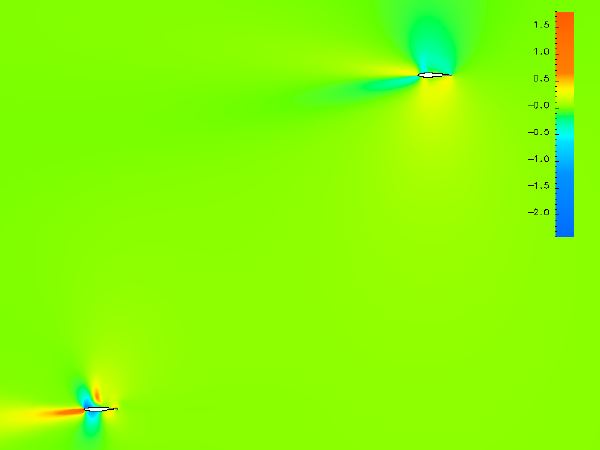}}
        \caption{$|\omega_l|:|\omega_d|= 1:10$}
    \end{subfigure}
    \\
    \begin{subfigure}{0.45\textwidth}
        \frame{\includegraphics[width=\textwidth]{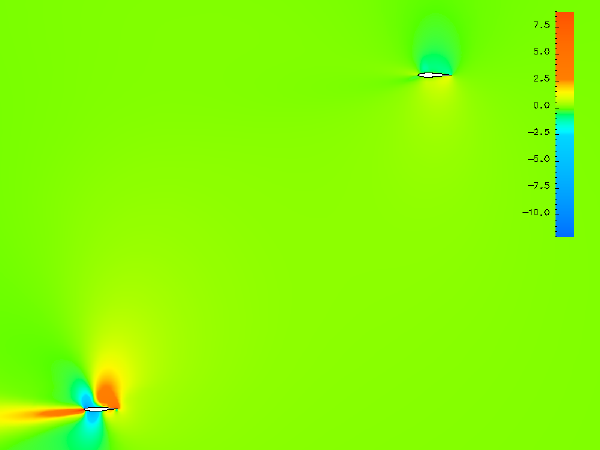}}
        \caption{$|\omega_l|:|\omega_d| = 1:50$}
    \end{subfigure}
    \hfill
    \begin{subfigure}{0.45\textwidth}
        \frame{\includegraphics[width=\textwidth]{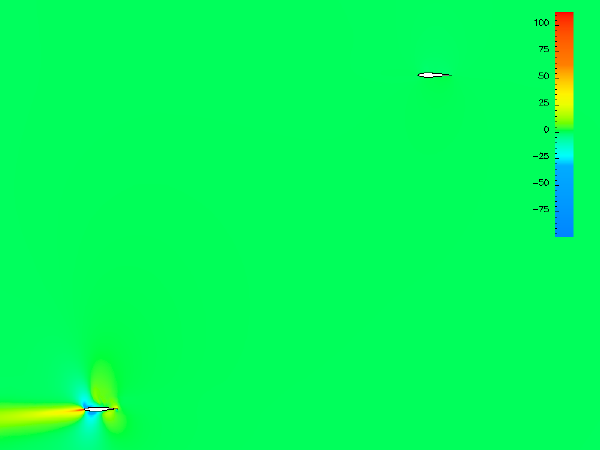}}
        \caption{$|\omega_l|:|\omega_d| = 1:500$}
    \end{subfigure}

    \caption{The dual solutions generated by the single mesh method where the coefficients of lift and drag are different.}
    \label{fourDual}
\end{figure}

The effect of different coefficients' ratios can be seen in Fig.\ref{fourDual}. The order of magnitude of lift is $43$ times larger than that of the drag. Thus, the dual solutions of lift will dominate the mesh adaptation process when $|\omega_l|:|\omega_d|= 1:1$. Then we try to increase the proportion of drag. While the dual solutions of drag will dominate the mesh adaptation if the coefficient is set too large, the suitable interval ranges from $|\omega_l|:|\omega_d| = 1:10$ to $|\omega_l|:|\omega_d|= 1:50$. As shown in Fig.\ref{isonLine}, the isolines of the single mesh method with coefficient $|\omega_l|:|\omega_d|= 1:50$ are smoother than other cases, which can be attributed to the ratio getting closer to the ratio of $F^1:F^2$. Meanwhile, the record of the single mesh method with coefficient $|\omega_l|:|\omega_d|= 1:50$ shown in Fig.\ref{compareSingle} behaves better than other cases for both lift and drag calculation. Nevertheless, the selection of coefficients introduces manual intervention. The coefficients based on the ratio of different $F^i$ are not known in advance for practical applications. Thus, it can not preserve accuracy during the shape optimization process where unexpected oscillations may destroy the adaptation profoundly. Besides, out of the motivation to conduct the shape optimization of airfoil, the dual solutions play an important role in ameliorating the shape that meets the optimization objectives. While the dual solutions generated from the single mesh method cannot be adopted to calculate useful shape derivatives that can improve the target functional, the optimization process may not reach the expected shape of the airfoil.

\subsection{Lift-drag ratio}
The lift-drag ratio is crucial for the shape optimal design of the airfoil. However, it is not easy to calculate correctly. The issues originate from the different orders of magnitude of lift and drag, the manual intervention of the algorithm, and the robustness during the mesh deformation process. It is shown in \cite{wang2024mechanisminformed} that the calculation of the lift-drag ratio may encounter data anomalies issues. With the multi-mesh method established above, we are able to develop a more robust algorithm for this multi-target functional. We conduct the algorithm with the following configurations.
\begin{itemize}
  \item  A domain surrounded by an outer
circle with a radius of 35;
\item  Mach number 0.8, and attack angle
1.25$^\circ$; 
\item Lax-Friedrichs numerical flux;
\item  NACA0012 at $(0,0)$;
\item Mirror
  reflection as the solid wall boundary condition;
\item Objectives: Lift-drag ratio of NACA0012.
\end{itemize}

\begin{figure}[!h]\centering
\includegraphics[width=0.6\textwidth]{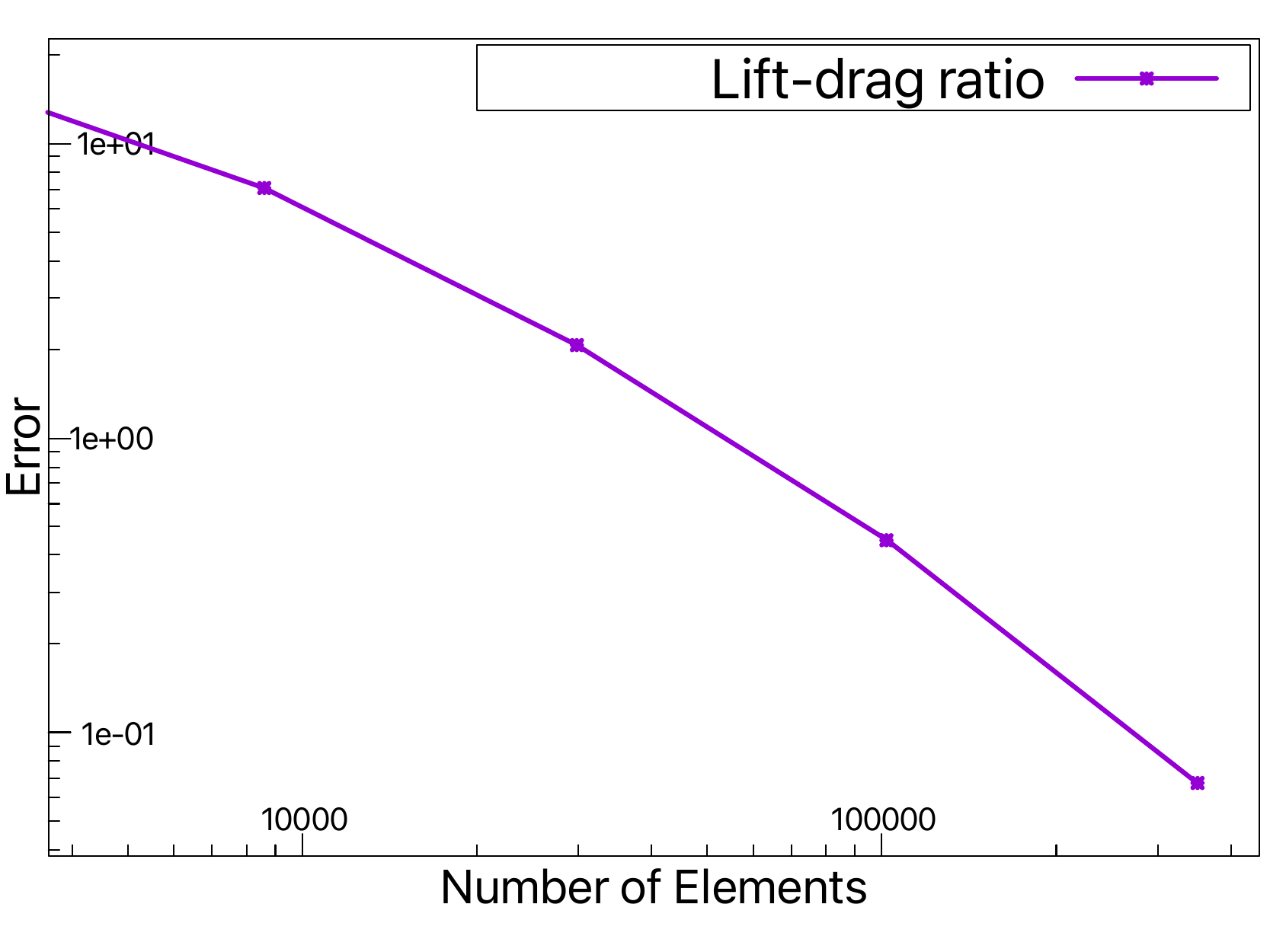}
\caption{The convergence of lift-drag ratio with the increase in the element size.}
\label{RatioRecord}
\end{figure}

\begin{figure}[!h]\centering
\frame{\includegraphics[width=0.48\textwidth]{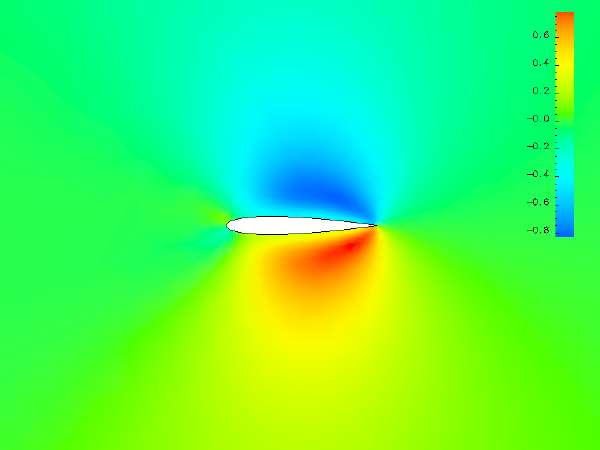}}  
\frame{\includegraphics[width=0.48\textwidth]{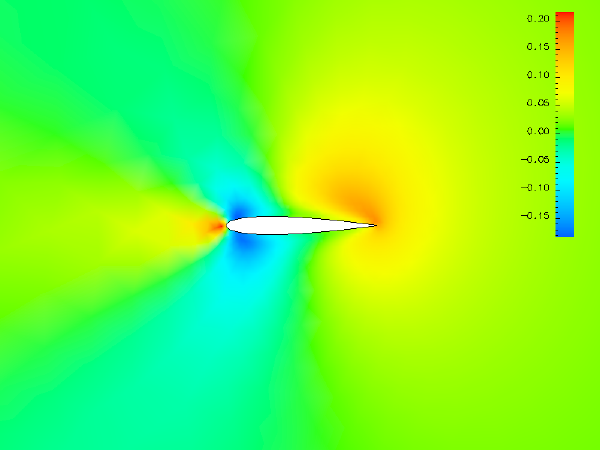}}
\caption{Left: dual solution of lift around NACA0012. Right: dual solution of drag around NACA0012.}
\label{NACA0012Ratio}
\end{figure}

Similarly, the initial mesh gets uniformly refined with $5$ rounds to produce a reference value, which is on a mesh with $3,723,264$ elements where lift is $0.406235$ while drag is $0.0259659$. The lift-drag ratio is $15.6449$. The dual solutions of both lift and drag are calculated on the same mesh as illustrated in Fig.\ref{NACA0012Ratio}. Moreover, a steady convergence of the lift-drag ratio can be obtained with the increase in the element size as shown in Fig. \ref{RatioRecord}.
In \cite{wang2024mechanisminformed}, the lift-drag ratio may not be robust if the mesh adaptation is conducted towards the drag. With the multi-mesh approach, a steady convergence can be guaranteed, which can mitigate the data abnormal issues significantly. As both drag and lift are derived in a dual-consistent manner, the multi-mesh version can preserve robustness for a complicated multi-target functional.

\begin{figure}[!h]\centering
\frame{\includegraphics[width=0.33\textwidth]{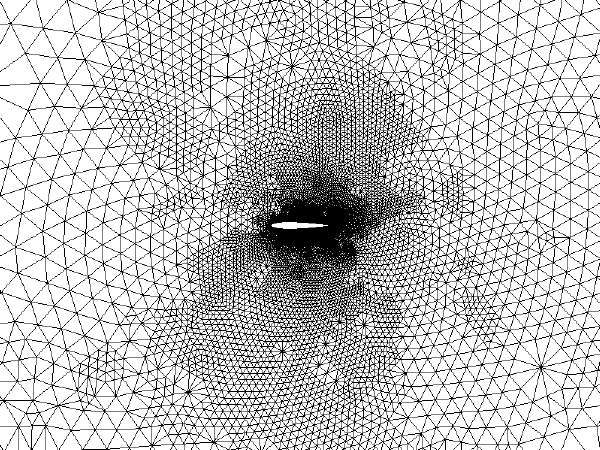}} 
\frame{\includegraphics[width=0.33\textwidth]{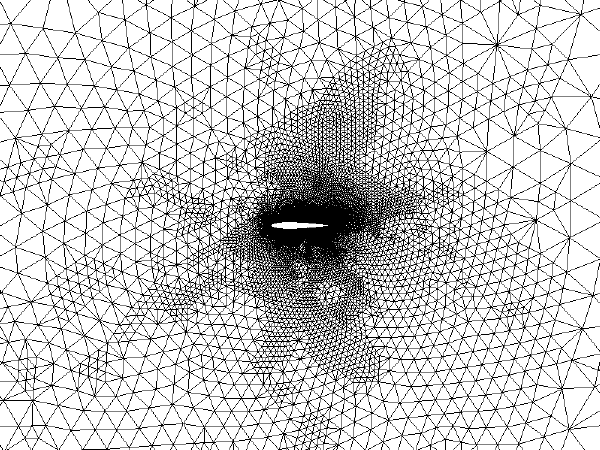}}
\frame{\includegraphics[width=0.33\textwidth]{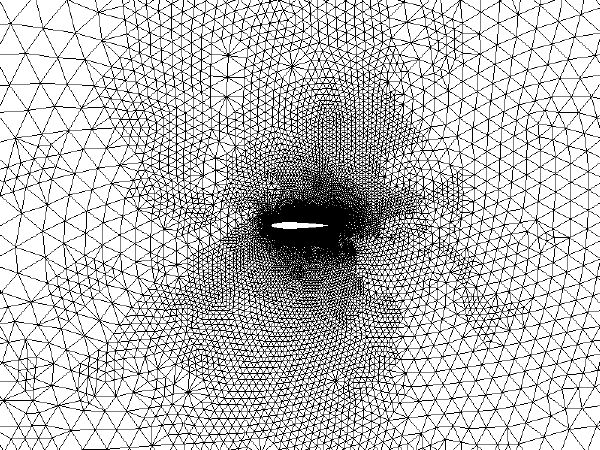}}
\caption{Mesh around the NACA0012. Left: generated from the lift. Middle: generated from drag. Right: common mesh.}
\label{LDprocess}
\end{figure} 

As shown in Fig. \ref{LDprocess}, the refined areas from the lift and drag are different even though they both demand the accuracy of the pressure alongside the airfoil. With the algorithm introduced above, we obtained the mesh for calculating the lift-drag ratio. Different from the single-mesh method where the dual equations are derived to take a balance between different components, this framework can preserve the accuracy of dual solutions concerning different target functionals. Such a framework provides a basic for further study on multi-target optimization. In order to demonstrate the efficiency of this multi-mesh framework, we conducted experiments in more complicated scenarios.

\subsection{Lift-Moment-Drag}
Besides lift and drag values, other values are also important for the aerodynamic design like the moment coefficient.
\begin{equation}\label{moment}
  \mathcal{J}_m(\mathbf{u})= \frac{1}{C_{\infty}}\int_{\Gamma}(\mathbf{x}-\mathbf{x}_{ref})\times p\mathbf{n}~ds,
\end{equation}
where $\mathbf{x}_{ref}$ is the moment reference point while $\mathbf{x}$ is the first coordinate of the integrated region along the airfoil.
In this part, we consider a scenario with three target functionals on three airfoils where the configurations are listed below.

\begin{itemize}
  \item  A domain surrounded by an outer
circle with a radius of 35;
\item  Mach number 0.8, and attack angle
1.25$^\circ$; 
\item Lax-Friedrichs numerical flux;
\item  NACA0012 at $(-2,0)$, RAE2822 at $(0,0)$, RAE2822 at $(2,0)$;
\item Mirror
  reflection as the solid wall boundary condition;
\item Objectives: Lift of the leading airfoil. Moment of the middle airfoil, $\mathbf{x}_{ref}=0.25$. Drag of the tail airfoil.
\end{itemize}

\begin{figure}[!h]\centering
\frame{\includegraphics[width=0.33\textwidth]{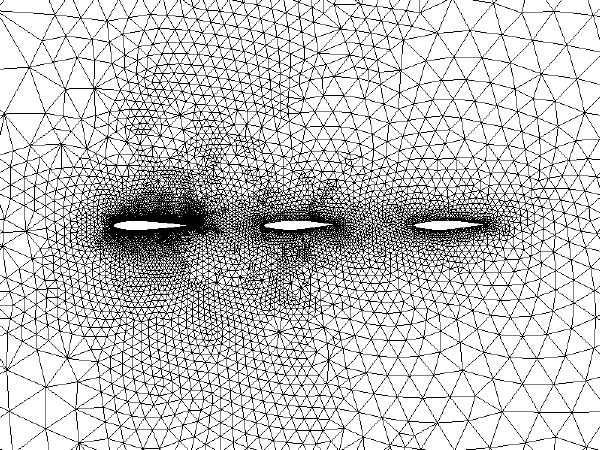}}
\frame{\includegraphics[width=0.33\textwidth]{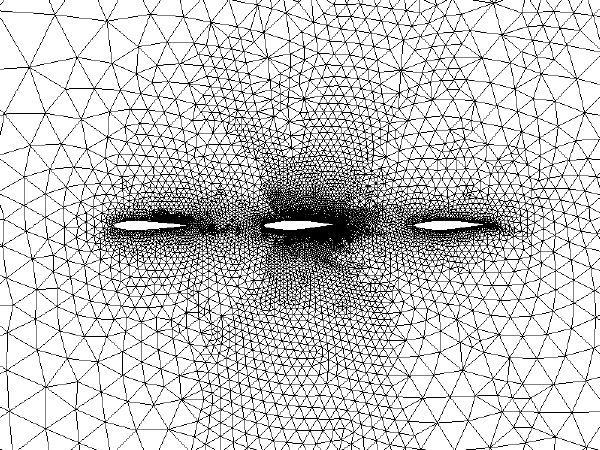}}
\frame{\includegraphics[width=0.33\textwidth]{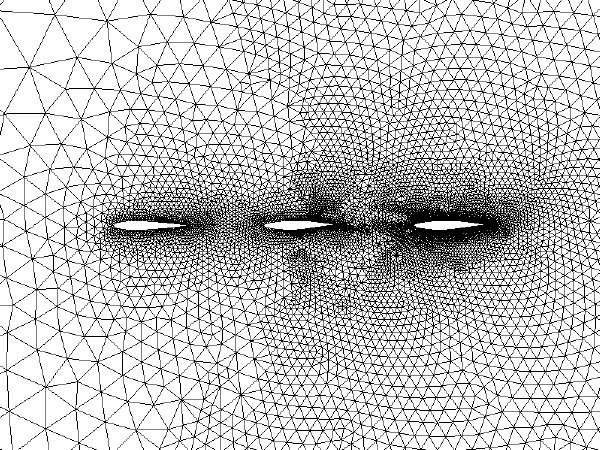}}
\caption{Meshes generated with different target functionals. Left: lift of the leading airfoil. Middle: the moment coefficient of the middle airfoil. Right: the drag of the tail airfoil.}
\label{ThreeMesh}
\end{figure}

Similarly, we obtain the reference value from the mesh with $5$ times uniform refinement, resulting in $9,379,840$ elements. The lift of the leading airfoil is $0.70737$, the moment coefficient of the middle airfoil is $0.17185$, and the drag of the tail airfoil is $0.013626$. In this three airfoils scenario, additional challenges occur from different aspects, i). the target functionals have different orders of magnitude, which makes the single-mesh method harder to select suitable coefficients. ii). the calculation of these three target functionals is not independent since the flows are influenced by other parts as well. The multi-mesh algorithm is extended to three target forms. Then after the first time DWR-based mesh adaptation for different target functional, three different meshses are generated as illustrated in Fig.\ref{ThreeMesh}.

\begin{figure}[!h]\centering
\frame{\includegraphics[width=0.45\textwidth]{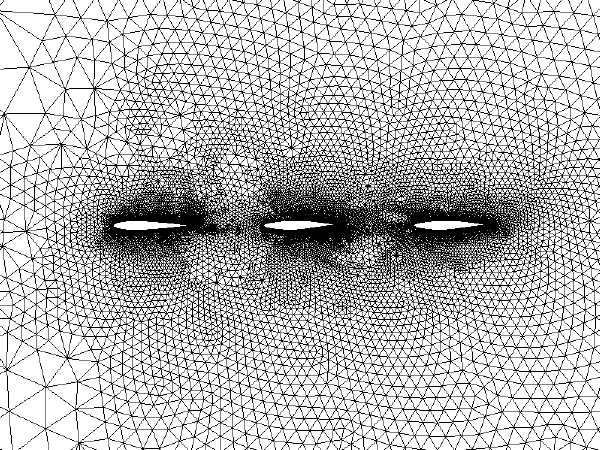}}
\includegraphics[width=0.45\textwidth]{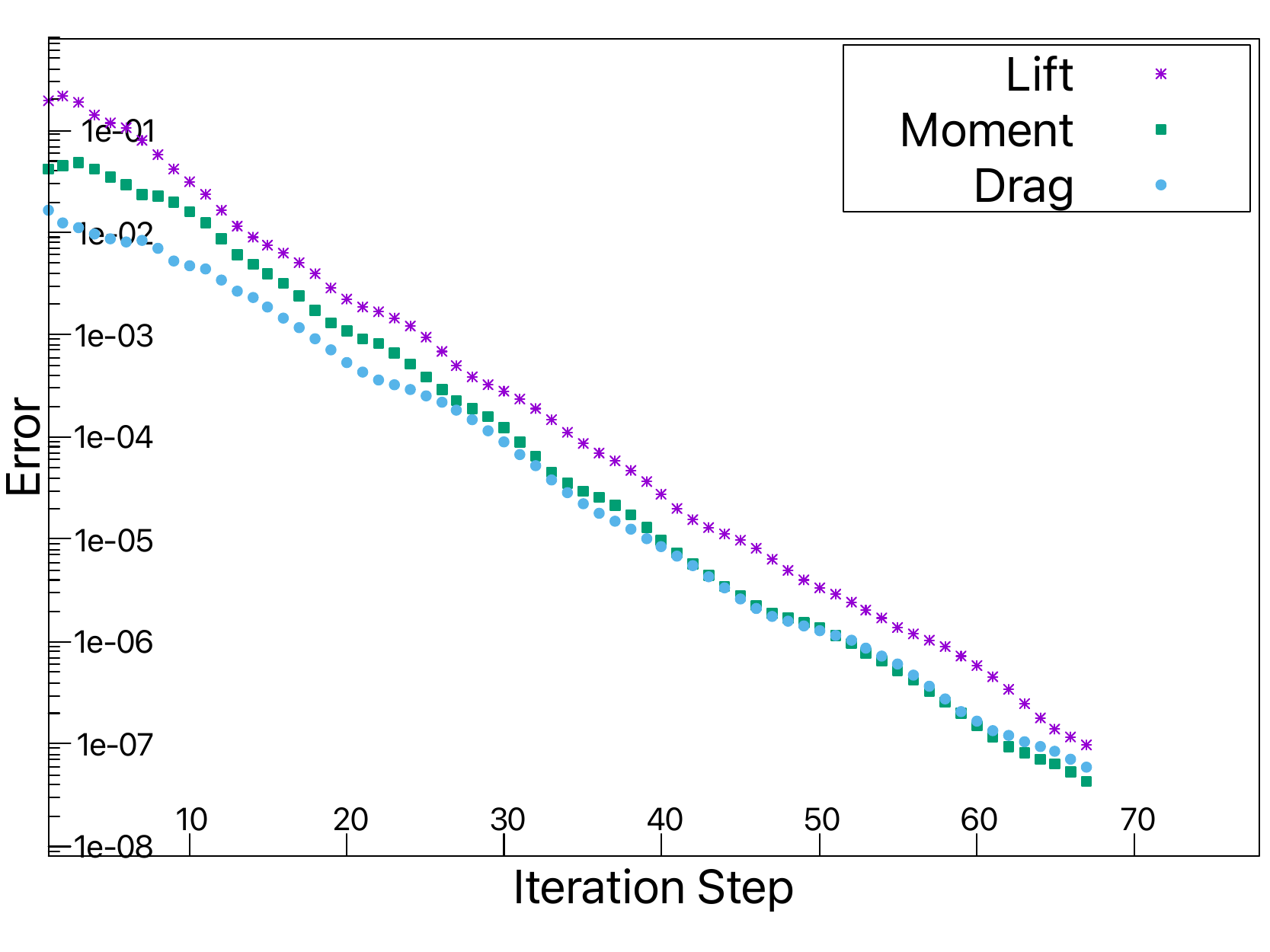}
\caption{Left: the common mesh generated with three different target functionals. Right: the iteration process of dual equations with different target functionals.}
\label{ThreeIter}
\end{figure}

Subsequently, the common mesh is generated with the hierarchical geometry tree as shown in Fig. \ref{ThreeIter}. The dual solutions are calculated simultaneously on this common mesh. The iterating history demonstrates that the three different dual equations can be solved well as their residuals' order of magnitude is comparable with the iteration processed.

\begin{figure}[!h]\centering
\frame{\includegraphics[width=0.33\textwidth]{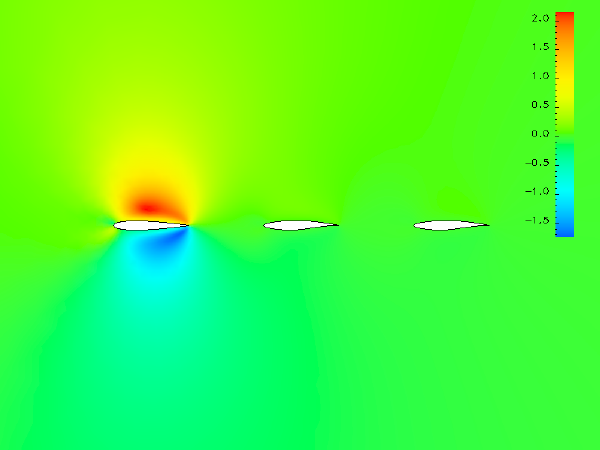}}
\frame{\includegraphics[width=0.33\textwidth]{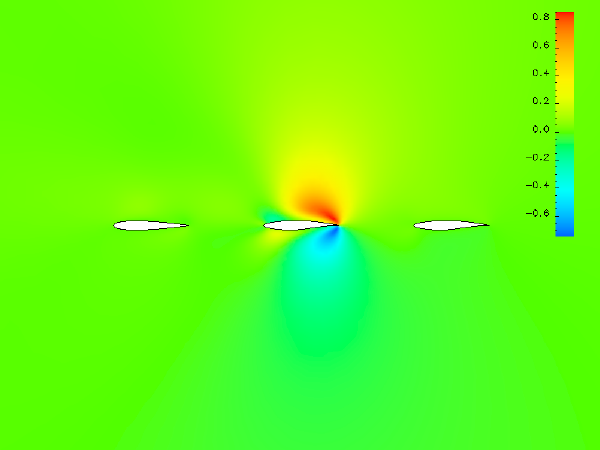}}
\frame{\includegraphics[width=0.33\textwidth]{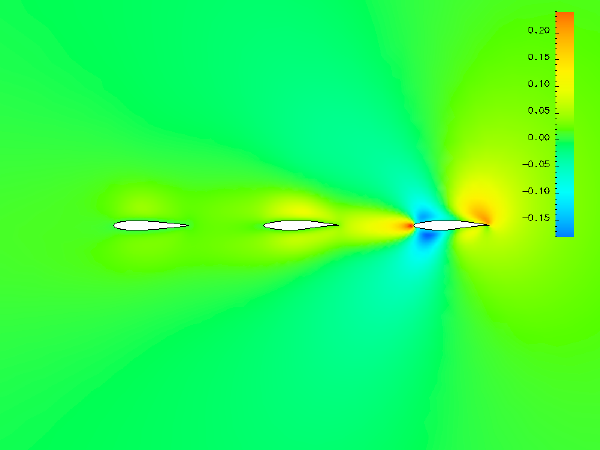}}
\caption{The distributions of the dual solutions concerning different target functionals. Left: lift of the leading airfoil. Middle: the moment coefficient of the middle airfoil. Right: the drag of the tail airfoil.}
\label{ThreeDual}
\end{figure} 

\begin{figure}[!h]\centering
\includegraphics[width=0.33\textwidth]{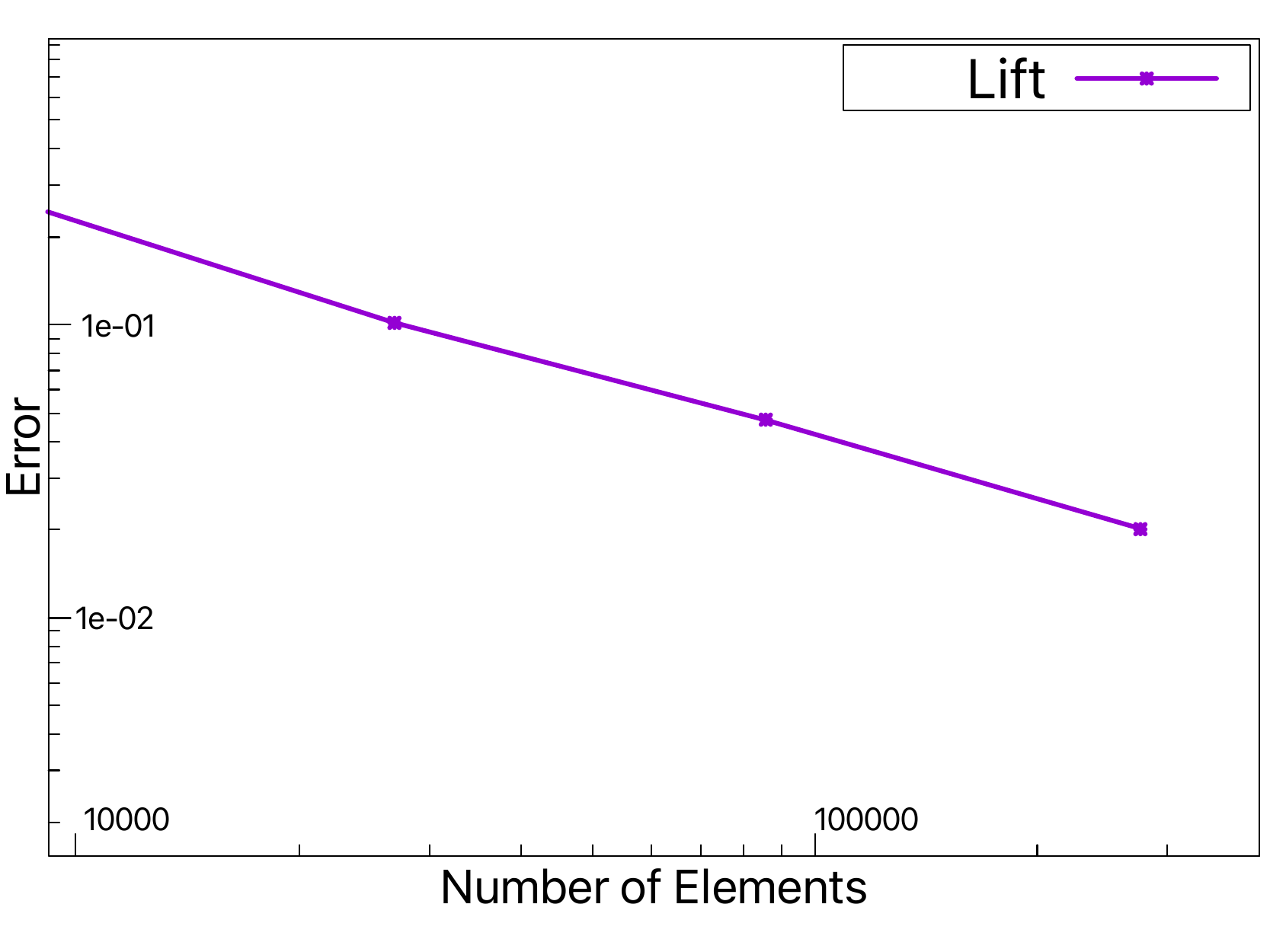}
\includegraphics[width=0.33\textwidth]{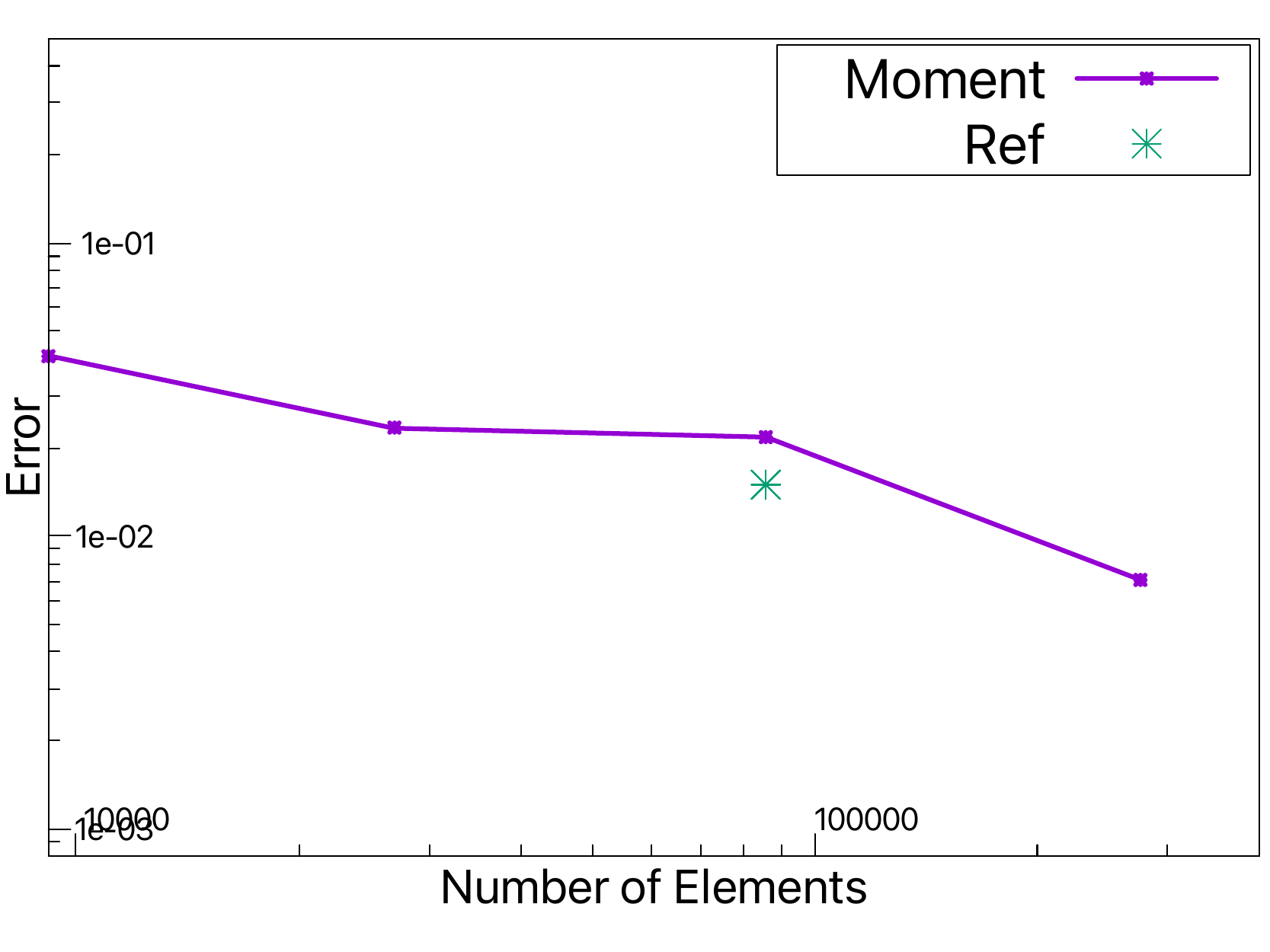}
\includegraphics[width=0.33\textwidth]{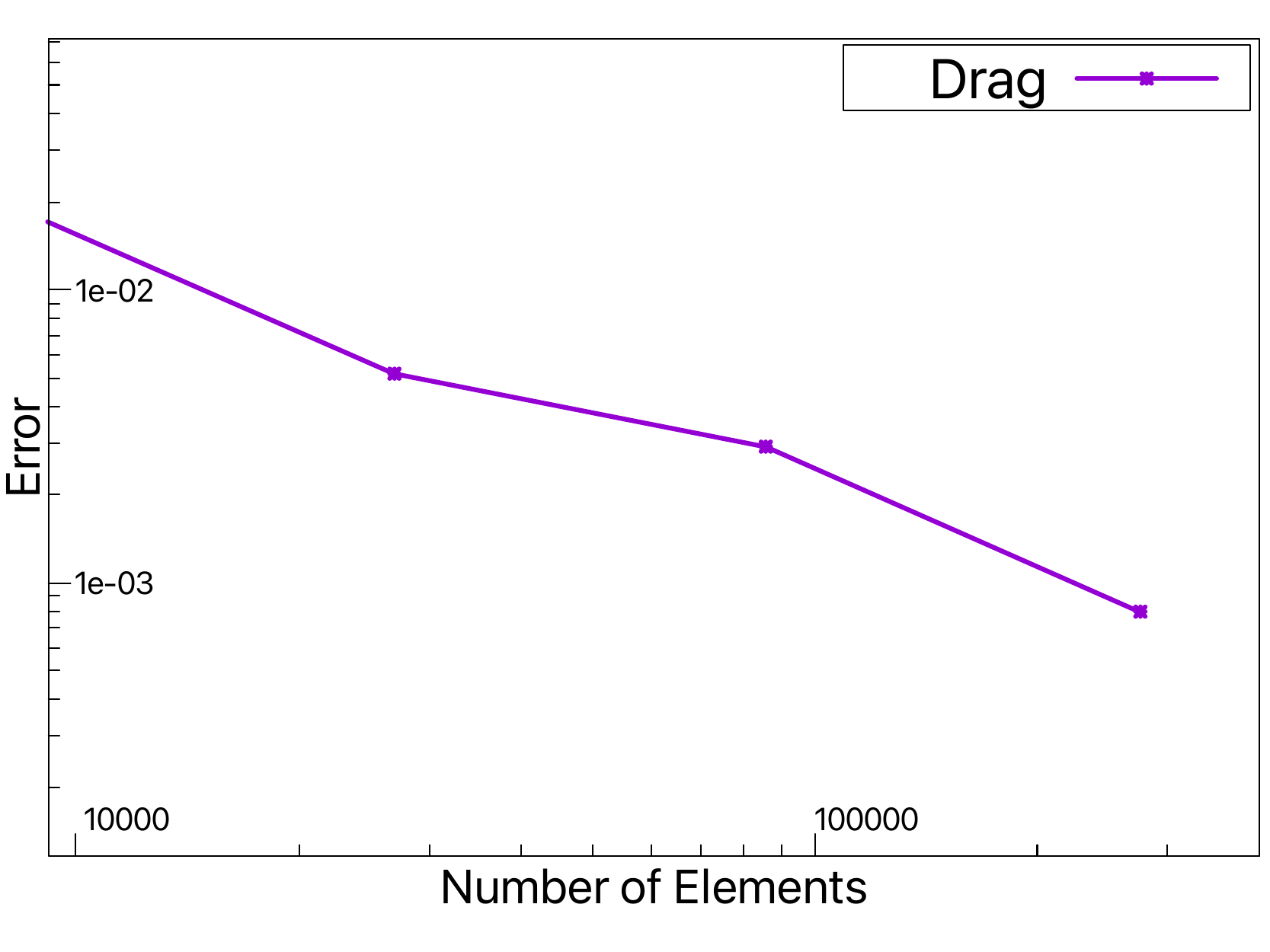}
\caption{The convergence of different target functionals. Left: lift of the leading airfoil. Middle: the moment coefficient of the middle airfoil. Right: the drag of the tail airfoil.}
\label{ThreeTarget}
\end{figure} 

While the dual solutions are based on different target functionals on different airfoils, they have different behaviors and corresponding dual solutions centered around different regions. As seen from Fig. \ref{ThreeDual}, the dual solutions of drag around the tail airfoil also demonstrate that the other two different airfoils need to be calculated well to generate a high quality mesh for calculating the target functional. With the multi-mesh DWR-based mesh refinement processed, we obtain finer meshes. From the result in Fig. \ref{ThreeTarget}, the error of target functional gets decreased when the element size gets increased. Even though the calculation of the moment of the middle airfoil exhibited oscillation during the refining process, the error still decreased when the refining level was stepped. The issues can be attributed to the middle airfoil is also influenced by the left two parts. Since the flow fields around the middle airfoil are not obtained precisely at the beginning, then the error of target functional is also influenced profoundly. Such an issue can be mitigated if the tolerance of the residual of governing equations is set lower. However, the multi-mesh DWR-based mesh adaptation method still helps the multi-target functional to converge robustly when the element size gets increased.

\section{Conclusion}
In this paper, we propose a novel multi-target DWR-based mesh adaptation method and adopt the multi-mesh approach to calculate the dual solutions of different target functionals. Such an algorithm can be easily implemented with the hierarchical geometry tree structure. From the numerical experiments, lift-drag ratio and other kinds of complicated multi-target functional can be obtained robustly. Meanwhile, the dual solutions can be iterated simultaneously without compromising their accuracy. With such an efficient algorithm, we plan to conduct multi-target optimization in the future. Besides, the multi-mesh algorithm has not been parallelized thoroughly. Different acceleration modules shall be implemented to enhance this framework in the future. Moreover, from the theoretical aspect, the extrapolation method in the original DWR-based mesh adaptation method can be extended to higher-order terms. Whether the convergence order can be guaranteed within the multi-target version is unknown yet. Even though the convergence of the lift-drag ratio is stable, further analysis shall be focused on the dual consistency of the lift-drag ratio as well.
 
\section*{Acknowledgments}
The work of Guanghui Hu is supported by The Science and Technology Development Fund, Macao SAR (No. 0082/2020/A2), National Natural Science Foundation of China (Nos. 11922120, 11871489), MYRG of University of Macau (No. MYRG-GRG2023-00157-FST-UMDF). The work of Ruo Li is partially supported by the National Natural Science Foundation of China (Grant No. 12288101). The work of Jingfeng Wang is supported by the UM Macao PhD Assistantship.
\bibliographystyle{plain}  
\bibliography{references}

\begin{thebibliography}{10}

\bibitem{arora2014lift}
Nipun Arora, Amit Gupta, Sanjeev Sanghi, Hikaru Aono, and Wei Shyy.
\newblock Lift-drag and flow structures associated with the “clap and
  fling” motion.
\newblock {\em Physics of Fluids}, 26(7), 2014.

\bibitem{bao2012h}
Gang Bao, Guanghui Hu, and Di~Liu.
\newblock An h-adaptive finite element solver for the calculations of the
  electronic structures.
\newblock {\em Journal of Computational Physics}, 231(14):4967--4979, 2012.

\bibitem{bouhlel2020scalable}
Mohamed~Amine Bouhlel, Sicheng He, and Joaquim~RRA Martins.
\newblock Scalable gradient--enhanced artificial neural networks for airfoil
  shape design in the subsonic and transonic regimes.
\newblock {\em Structural and Multidisciplinary Optimization},
  61(4):1363--1376, 2020.

\bibitem{AFEPack}
{Z}henning {C}ai, Yun {C}hen, {Y}ana {D}i, {G}uanghui {H}u, {R}uo {L}i,
  {W}enbin {L}iu, {H}eyu {W}ang, {F}anyi {Y}ang, {C}hengbao {Y}ao, and
  {H}ongfei {Z}han.
\newblock {A}{F}{E}{P}ack: a general-purpose {C}++ library for numerical
  solutions of partial differential equations.
\newblock {\em Communications in Computational Physics}, In press, 2024.

\bibitem{cary2021cfd}
A~Cary, J~Chawner, E~Duque, W~Gropp, B~Kleb, R~Kolonay, E~Nielsen, and B~Smith.
\newblock {The CFD vision 2030 roadmap: 2020 status progress and challenges}.
\newblock {\em AIAA Paper}, 2726:2021, 2021.

\bibitem{cary2022realizing}
Andrew Cary, John Chawner, Earl Duque, William Gropp, Bil Kleb, Ray Kolonay,
  Eric Nielsen, and Brian Smith.
\newblock Realizing the vision of {CFD} in 2030.
\newblock {\em Computing in Science \& Engineering}, 24(1):64--70, 2022.

\bibitem{chen2020output}
{G}uodong {C}hen and {K}rzysztof {F}idkowski.
\newblock {O}utput-based error estimation and mesh adaptation using
  convolutional neural networks: {A}pplication to a scalar advection-diffusion
  problem.
\newblock In {\em {A}{I}{A}{A} {S}citech 2020 {F}orum}, page 1143, 2020.

\bibitem{chen2021output}
{G}uodong {C}hen and {K}rzysztof~{J} {F}idkowski.
\newblock {Ou}tput-based adaptive aerodynamic simulations using convolutional
  neural networks.
\newblock {\em {C}omputers \& {F}luids}, 223:104947, 2021.

\bibitem{di2009computation}
Yana Di and Ruo Li.
\newblock Computation of dendritic growth with level set model using a
  multi-mesh adaptive finite element method.
\newblock {\em Journal of Scientific Computing}, 39:441--453, 2009.

\bibitem{DOLEJSI2021178}
V{\'\i}t Dolej{\v{s}}{\'\i}, Ond{\v{r}}ej Barto{\v{s}}, and Filip Roskovec.
\newblock Goal-oriented mesh adaptation method for nonlinear problems including
  algebraic errors.
\newblock {\em Computers \& Mathematics with Applications}, 93:178--198, 2021.

\bibitem{dolejvsi2023anisotropic}
V{\'\i}t Dolej{\v{s}}{\'\i} and Georg May.
\newblock An anisotropic hp-mesh adaptation method for time-dependent problems
  based on interpolation error control.
\newblock {\em Journal of Scientific Computing}, 95(2):36, 2023.

\bibitem{dolejsi2022}
Vít Dolejší and Filip Roskovec.
\newblock Goal-oriented anisotropic hp-adaptive {D}iscontinuous {G}alerkin
  {M}ethod for the {E}uler {E}quations.
\newblock {\em Communications on Applied Mathematics and Computation},
  4:143--179, 2022.

\bibitem{endtmayer2019mesh}
Bernhard Endtmayer, Ulrich Langer, Ira Neitzel, Thomas Wick, and Winnifried
  Wollner.
\newblock {Mesh adaptivity and error estimates applied to a regularized
  p-Laplacian constrainted optimal control problem for multiple quantities of
  interest}.
\newblock {\em PAMM}, 19(1):e201900231, 2019.

\bibitem{endtmayer2024posteriori}
Bernhard Endtmayer, Ulrich Langer, Thomas Richter, Andreas Schafelner, and
  Thomas Wick.
\newblock A posteriori single-and multi-goal error control and adaptivity for
  partial differential equations.
\newblock {\em ArXiv preprint arXiv:2404.01738}, 2024.

\bibitem{gramacy2020surrogates}
Robert~B Gramacy.
\newblock {\em Surrogates: Gaussian process modeling, design, and optimization
  for the applied sciences}.
\newblock Chapman and Hall/CRC, 2020.

\bibitem{hartmann2007adjoint}
Ralf Hartmann.
\newblock {A}djoint consistency analysis of discontinuous {G}alerkin
  discretizations.
\newblock {\em SIAM Journal on Numerical Analysis}, 45(6):2671--2696, 2007.

\bibitem{hartmann2008multitarget}
Ralf Hartmann.
\newblock Multitarget error estimation and adaptivity in aerodynamic flow
  simulations.
\newblock {\em SIAM Journal on Scientific Computing}, 31(1):708--731, 2008.

\bibitem{hartmann2002adaptive}
Ralf Hartmann and Paul Houston.
\newblock Adaptive discontinuous {G}alerkin finite element methods for the
  compressible {E}uler equations.
\newblock {\em Journal of Computational Physics}, 183(2):508--532, 2002.

\bibitem{HARTMANN2015754}
Ralf Hartmann and Tobias Leicht.
\newblock Generalized adjoint consistent treatment of wall boundary conditions
  for compressible flows.
\newblock {\em Journal of Computational Physics}, 300:754--778, 2015.

\bibitem{hu2013adaptive}
Guanghui Hu.
\newblock An adaptive finite volume method for 2{D} steady {E}uler equations
  with {WENO} reconstruction.
\newblock {\em Journal of Computational Physics}, 252:591--605, 2013.

\bibitem{hu2010jcp}
Guanghui Hu, Ruo Li, and Tao Tang.
\newblock A robust high-order residual distribution type scheme for steady
  {E}uler equations on unstructured grids.
\newblock {\em Journal of Computational Physics}, 229(5):1681--1697, 2010.

\bibitem{hu2011robust}
Guanghui Hu, Ruo Li, and Tao Tang.
\newblock A robust {WENO} type finite volume solver for steady {E}uler
  equations on unstructured grids.
\newblock {\em Communications in Computational Physics}, 9(3):627--648, 2011.

\bibitem{hu2016adjoint}
Guanghui Hu, Xucheng Meng, and Nianyu Yi.
\newblock Adjoint-based an adaptive finite volume method for steady {E}uler
  equations with non-oscillatory k-exact reconstruction.
\newblock {\em Computers \& Fluids}, 139:174--183, 2016.

\bibitem{HU2016235}
Guanghui Hu and Nianyu Yi.
\newblock An adaptive finite volume solver for steady {E}uler equations with
  non-oscillatory k-exact reconstruction.
\newblock {\em Journal of Computational Physics}, 312:235--251, 2016.

\bibitem{hu2009multi}
Xianliang Hu, Ruo Li, and Tao Tang.
\newblock A multi-mesh adaptive finite element approximation to phase field
  models.
\newblock {\em Communications in Computational Physics}, 5(5):1012--1029, 2009.

\bibitem{jiang2020surrogate}
Ping Jiang, Qi~Zhou, Xinyu Shao, Ping Jiang, Qi~Zhou, and Xinyu Shao.
\newblock {\em Surrogate-model-based design and optimization}.
\newblock Springer, 2020.

\bibitem{kuang2023mul}
Yang Kuang, Yedan Shen, and Guanghui Hu.
\newblock Towards chemical accuracy using a multi-mesh adaptive finite element
  method in all-electron density functional theory.
\newblock {\em ArXiv preprint arXiv:2310.15651}, 2023.

\bibitem{lee2020wind}
Sang-Lae Lee and SangJoon Shin.
\newblock Wind turbine blade optimal design considering multi-parameters and
  response surface method.
\newblock {\em Energies}, 13(7):1639, 2020.

\bibitem{leifsson2016multiobjective}
Leifur Leifsson, Slawomir Koziel, and Yonatan~A Tesfahunegn.
\newblock Multiobjective aerodynamic optimization by variable-fidelity models
  and response surface surrogates.
\newblock {\em AIAA Journal}, 54(2):531--541, 2016.

\bibitem{li2005multi}
Ruo Li.
\newblock On multi-mesh h-adaptive methods.
\newblock {\em Journal of Scientific Computing}, 24:321--341, 2005.

\bibitem{li2008multigrid}
Ruo Li, Xin Wang, and Weibo Zhao.
\newblock A multigrid block {LU-SGS} algorithm for {E}uler equations on
  unstructured grids.
\newblock {\em Numerical Mathematics: Theory, Methods and Applications},
  1:92--112, 2008.

\bibitem{liu2022mp}
Chengyu Liu and Guanghui Hu.
\newblock {A}n {M}{P}-{D}{W}{R} method for h-adaptive finite element methods.
\newblock {\em Numerical Algorithms}, pages 1--21, 2023.

\bibitem{martins2022aerodynamic}
Joaquim~RRA Martins.
\newblock Aerodynamic design optimization: Challenges and perspectives.
\newblock {\em Computers \& Fluids}, 239:105391, 2022.

\bibitem{meng2021fourth}
Xucheng Meng, Yaguang Gu, and Guanghui Hu.
\newblock A fourth-order unstructured {NURBS}-enhanced finite volume {WENO}
  scheme for steady {E}uler equations in curved geometries.
\newblock {\em Communications on Applied Mathematics and Computation}, pages
  1--28, 2021.

\bibitem{Meng2022}
Xucheng Meng and Guanghui Hu.
\newblock A nurbs-enhanced finite volume method for steady euler equations with
  goal-oriented $h$-adaptivity.
\newblock {\em Communications in Computational Physics}, 32:490--523, 06 2022.

\bibitem{nemec2014toward}
{M}arian {N}emec and {M}ichael {A}ftosmis.
\newblock {T}oward automatic verification of goal-oriented flow simulations.
\newblock Technical Report {N}{A}{S}{A}/{T}{M}-2014-218386, 2014.

\bibitem{ojha2021adaptive}
Vivek Ojha, Krzysztof Fidkowski, and Carlos~E Cesnik.
\newblock Adaptive mesh refinement for fluid-structure interaction simulations.
\newblock In {\em AIAA SciTech 2021 Forum}, page 0731, 2021.

\bibitem{ojstersek2020multi}
Robert Ojstersek, Miran Brezocnik, and Borut Buchmeister.
\newblock Multi-objective optimization of production scheduling with
  evolutionary computation: A review.
\newblock {\em International Journal of Industrial Engineering Computations},
  11(3):359--376, 2020.

\bibitem{pardo2010multigoal}
David Pardo.
\newblock Multigoal-oriented adaptivity for hp-finite element methods.
\newblock {\em Procedia Computer Science}, 1(1):1953--1961, 2010.

\bibitem{sforza2020direct}
Pasquale~M Sforza.
\newblock Direct calculation of zero-lift drag coefficients and ({L/D}) max in
  subsonic cruise.
\newblock {\em Journal of Aircraft}, 57(6):1224--1228, 2020.

\bibitem{singh2024hybrid}
Gyan Singh and Amit~K Chaturvedi.
\newblock Hybrid modified particle swarm optimization with genetic algorithm
  ({GA}) based workflow scheduling in cloud-fog environment for multi-objective
  optimization.
\newblock {\em Cluster Computing}, 27(2):1947--1964, 2024.

\bibitem{slotnick2014cfd}
Jeffrey~P Slotnick, Abdollah Khodadoust, Juan Alonso, David Darmofal, William
  Gropp, Elizabeth Lurie, and Dimitri~J Mavriplis.
\newblock {CFD} vision 2030 study: a path to revolutionary computational
  aerosciences.
\newblock Technical Report {NASA/CR}–2014-218178, 2014.

\bibitem{solin2010monolithic}
Pavel Sol{\'\i}n, J~Cerveny, Lenka Dubcova, and David Andrs.
\newblock Monolithic discretization of linear thermoelasticity problems via
  adaptive multimesh hp-{FEM}.
\newblock {\em Journal of Computational and Applied Mathematics},
  234(7):2350--2357, 2010.

\bibitem{vasile2020aerodynamic}
Joseph~D Vasile, Joshua Bryson, and Frank Fresconi.
\newblock Aerodynamic design optimization of long range projectiles using
  missile {DATCOM}.
\newblock In {\em AIAA Scitech 2020 Forum}, page 1762, 2020.

\bibitem{venditti2000adjoint}
David Venditti and David Darmofal.
\newblock Adjoint error estimation and grid adaptation for functional outputs:
  {A}pplication to quasi-one-dimensional flow.
\newblock {\em Journal of Computational Physics}, 164(1):204--227, 2000.

\bibitem{wang2024towards}
Jingfeng Wang and Guanghui Hu.
\newblock {Towards the efficient calculation of quantity of interest from
  steady Euler equations II: a CNNs-based automatic implementation}.
\newblock {\em Communications in Computational Physics}, In press, ArXiv
  preprint:2308.07140, 2023.

\bibitem{wang2024mechanisminformed}
Jingfeng Wang and Guanghui Hu.
\newblock A mechanism-driven reinforcement learning framework for shape
  optimization of airfoils.
\newblock {\em ArXiv preprint}, 2403.04329, 2024.

\bibitem{wang2023towards}
{J}ingfeng {W}ang and {G}uanghui {H}u.
\newblock {T}owards the efficient calculation of quantity of interest from
  steady {E}uler equations {I}: a dual-consistent {D}{W}{R}-based h-adaptive
  {N}ewton-{G}{M}{G} solver.
\newblock {\em Communications in Computational Physics}, 35(3):579--608, 2024.

\bibitem{wu2022gradient}
Neil Wu, Charles~A Mader, and Joaquim~RRA Martins.
\newblock A gradient-based sequential multifidelity approach to
  multidisciplinary design optimization.
\newblock {\em Structural and Multidisciplinary Optimization}, 65(4):131, 2022.

\bibitem{wu2018multi}
Shengyang Wu, Xianliang Hu, and Shengfeng Zhu.
\newblock A multi-mesh finite element method for phase-field based photonic
  band structure optimization.
\newblock {\em Journal of Computational Physics}, 357:324--337, 2018.

\end{thebibliography}

\end{document}